\documentclass[12pt]{article}

\title
{Passage of L\'{e}vy Processes across Power Law Boundaries at
Small 
Times\footnote{This research was partially supported by ARC
grants
DP0210572 and DP0664603}}
\author{J. Bertoin
\thanks{Laboratoire des Probabilit\'{e}s, Universit\'{e}
 Pierre et Marie Curie,
Paris, France; email:jbe@ccr.jussieu.fr}
\and
R. A. Doney
\thanks{Department of Mathematics, University of Manchester,
Manchester M13 9PL, UK; email:rad@ma.man.ac.uk.}
\and
R. A. Maller
\thanks{Centre for Mathematical Analysis, and School of Finance
and Applied Statistics, Australian National University,
Canberra, ACT;
email:Ross.Maller@anu.edu.au
}}
\date{}

\usepackage{amssymb}
\usepackage{amscd}
\usepackage{amsmath}%
\usepackage{amsfonts}%
\usepackage{graphicx}

\newtheorem{remark}{Remark}

\numberwithin{equation}{section}
\newtheorem{thm}{Theorem}[section]
\newtheorem{lem}{Lemma}[section]
\newtheorem{prop}{Proposition}[section]
\newtheorem{cor}{Corollary}[section]

\def\todr{\buildrel D \over \to }

\def\topr{\buildrel P \over \to }
\def\veps{\varepsilon}
\def\bbox{\hfil \vrule height5pt width4pt depth0pt}
\def\hat{\widehat}
\def\tilde{\widetilde}
\def\pibar{\overline{\Pi}}

\def\R{\mathbb{R}}

\def\E{\mathbb{E}}
\def\P{\mathbb{P}}

\def\D{\mathbb{D}}

\def \e{{\rm e}}

\begin{document}
\maketitle

\begin{abstract}
We wish to characterise when a   L\'{e}vy process $X_t$ crosses
boundaries like $t^\kappa$, $\kappa>0$,
in a one or two-sided sense, for
small times $t$; thus, we enquire when
$\limsup_{t\downarrow 0}|X_t|/t^{\kappa}$,
$\limsup_{t\downarrow 0}X_t/t^{\kappa}$ 
and/or
$\liminf_{t\downarrow 0}X_t/t^{\kappa}$ 
are almost surely (a.s.) finite or infinite.
Necessary and sufficient 
conditions are given for  these possibilities
for all values of $\kappa>0$. Often (for many values of
$\kappa$), when the limsups are finite
a.s., they are in fact zero,  as we show,
but the limsups may in some circumstances take finite, 
nonzero, values, a.s.
In general, the process crosses one or two-sided boundaries in
quite different ways, but surprisingly this is not so for 
the case $\kappa=1/2$.
An integral test is given to distinguish the possibilities
in that case.
Some results relating to other norming sequences for $X$, and
when $X$ is centered at a nonstochastic function, are also given.
\end{abstract}

\vfill

\noindent
\begin{tabbing}
\\
\> secondary:  60G10, 60G17, 60J65.
\end{tabbing}
\vspace{1cm}


\newpage

\section{Introduction}\label{s1}\
Let $X=(X_{t},t\geq0)$ be a L\'{e}vy process with
characteristic triplet $(\gamma,\sigma,\Pi)$,
where $\gamma\in \R$, $\sigma^2\ge 0$, and the L\'evy
measure $\Pi$ has $\int(x^2\wedge 1) \Pi({\rm d}x)$ finite.
See \cite{b} and \cite{satoa} for basic definitions and
properties.   
Since we will only be concerned with local behaviour of $X_t$,
as
$t\downarrow 0$, we can ignore the ``big jumps" in $X$ (those
with modulus exceeding 1, say),
  and write its characteristic exponent,
$\Psi(\theta) =\frac{1}{t}\log E\e^{{\rm i}\theta X_t}$, $\theta
\in\R$, as
\begin{equation}\label{canexp}
\Psi(\theta) ={\rm i}\gamma\theta -\sigma^2\theta^2/2 
+\int_{[-1,1]} (\e^{{\rm i}\theta x} -1-{\rm i}\theta x)
\Pi({\rm d}x)
. \end{equation}
The L\'evy process is of bounded variation, for which we use the
notation $X\in bv$, if and only if $\sigma^2=0$ and 
$\int_{|x|\le 1} |x|\Pi({\rm d}x)<\infty$, and in that case, we
denote by
$$\delta:=\gamma - \int_{[-1,1]}x\Pi({\rm d}x)$$
its drift coefficient.

We will continue the work of Blumenthal and Getoor \cite{bg} and
Pruitt \cite{pruitt}, in a sense, by investigating the possible
limiting
values taken by $t^{-\kappa}X_t$ as $t\downarrow 0$, where
$\kappa>0$
is some parameter.
Recall that 
Blumenthal and Getoor \cite{bg} introduced the upper-index
$$\beta:=\inf\left\{\alpha>0: \int_{|x|\leq
1}|x|^{\alpha}\Pi({\rm d}x)<\infty\right\}\in[0,2],$$
which plays a critical role in this framework. Indeed, assuming
for simplicity that
the Brownian coefficient $\sigma^2$ is zero, and further that
the drift coefficient
$\delta$ is also 0
when $X\in bv$, then with probability one, 
\begin{equation}\label{eqbg}
\limsup_{t\downarrow 0} \frac{|X_{t}|}{t^{\kappa}}=\left\{
\begin{matrix}
0\\
\infty\\
\end{matrix}\right. 
\hbox{ according as } \left\{ \begin{matrix}\kappa <1/\beta&\\
\kappa>1/\beta.&\\
\end{matrix}\right.
\end{equation}
See also  Pruitt \cite{pruitt}. Note that the critical case when
$\kappa=1/\beta$
is not covered by (\ref{eqbg})

One application of this kind of study is to get information 
on the rate of growth  of the process relative to
power law functions, in both a one and a two-sided sense, at
small times. More precisely, we are concerned with the values of
$\limsup_{t\downarrow 0}|X_t|/t^{\kappa}$ and of
$\limsup_{t\downarrow 0}X_t/t^{\kappa}$
(and the behaviour of $\liminf_{t\downarrow 0}X_t/t^{\kappa}$
can be deduced from the limsup behaviour by performing a sign
reversal).
 For example, when
\begin{equation}\label{main}
\limsup_{t\downarrow 0}\frac{|X_t|}{t^\kappa}=\infty\ {\rm
a.s.},
\end{equation}
then the regions
$\{(t,y)\in [0,\infty)\times\R: |y|> at^\kappa\}$ 
are entered infinitely often for arbitrarily small $t$, a.s.,
for
all $a>0$. This can be thought of as a kind of ``regularity" of
$X$ for these regions, at 0.
We will refer to this kind of behaviour 
as crossing a ``two-sided" boundary.
On the other hand, when
\begin{equation}\label{main2}
\limsup_{t\downarrow 0}\frac{X_t}{t^\kappa}=\infty\ {\rm a.s.},
\end{equation}
we have ``one-sided" (up)crossings;
and similarly for downcrossings, phrased in terms of the liminf.
In general, the process crosses one or two-sided boundaries in
quite different ways,
and, often (for many values of
$\kappa$), when the limsups are finite
a.s., they are in fact zero, a.s.,  as we will show.
But the limsups may in some circumstances take finite, 
nonzero, values, a.s.
Our aim here is to give necessary and sufficient 
conditions (NASC) which distinguish all these possibilities,
for all values of $\kappa>0$. 

Let us eliminate at the outset certain cases which are trivial
or easily deduced from known results.
A result of Khintchine \cite{khin}
(see Sato (\cite{satoa}, Prop. 47.11, p. 358))
is that, for any L\'evy process $X$ with Brownian coefficient
$\sigma^2\ge 0$, we have
\begin{equation}\label{khin}
\limsup_{t\downarrow 0}\frac{|X_t|}
{\sqrt{2t\log|\log t|}}=\sigma\ {\rm a.s.}
\end{equation}
Thus we immediately see that 
\eqref{main}  and \eqref{main2} cannot hold for $0< \kappa<1/2$;
we always have
$\lim_{t\downarrow 0}X_t/t^\kappa=0$ a.s. in these cases.
Of course, this also agrees with (\ref{eqbg}), since, always,
$\beta\le 2$.
More precisely, recall the decomposition $X_t=\sigma B_t+
X^{(0)}_t$,
where $X^{(0)}$ is a L\'evy process with characteristics
$(\gamma, 0, \Pi)$
and $B$ is an independent Brownian motion. Khintchine's law of the
iterated logarithm
for $B$ and \eqref{khin} applied for $X^{(0)}$ give
\begin{equation}\label{khin2}
-\liminf_{t\downarrow 0}\frac{X_t} {\sqrt{2t\log|\log t|}}=
\limsup_{t\downarrow 0}\frac{X_t} {\sqrt{2t\log|\log
t|}}=\sigma\ 
{\rm a.s.}
\end{equation}
So the one and 
two-sided limsup behaviours of $X$ are precisely
determined when $\sigma^2>0$ (regardless 
of the behaviour of $\Pi(\cdot)$, which may not even be
present).
With these considerations, it is clear that
{\it throughout, we can assume}
\begin{equation}\label{nohold}
\sigma^2=0\,.
 \end{equation}
{\it Furthermore, we can restrict attention to the cases
$\kappa\ge 1/2$.}

 A
result of Shtatland \cite{shtat} and Rogozin \cite{rog} is that
 $X\notin bv$ if and only if
\[ 
-\liminf_{t\downarrow 0}\frac{X_t}{t}
= \limsup_{t\downarrow 0}\frac{X_t}{t}=\infty\ {\rm  a.s.}, 
\]
so \eqref{main} and \eqref{main2} hold for all $\kappa\ge 1$,
in this case (and similarly for the liminf).
On the other hand, when $X\in bv$, we have
\[ \lim_{t\downarrow 0}\frac{X_t}{t}  =\delta, \  {\rm  a.s.},
\]
where 
$\delta$ is the drift of $X$ (cf. \cite{b}, p.84).
Thus if $\delta>0$, \eqref{main2} holds for all $\kappa> 1$,
but for no $\kappa\le 1$,
while if $\delta<0$, 
\eqref{main2} can  hold for no $\kappa>0$;
while \eqref{main} holds in either case, with $\kappa>1$,
but for no $\kappa\le 1$.
{\it Thus, when $X\in bv$, we need only consider the case
$\delta=0$.}

The main statements for two-sided (respectively, one-sided)
boundary crossing
will be given in Section \ref{s2} (respectively, Section
\ref{s3}) and proved in Section \ref{s4} (respectively,
Section \ref{s5}).
We use throughout similar notation to
\cite{dm1},  \cite{dm2} and  \cite{dm3}. In particular, we write
 $\Pi^{\#}$ for the L\'evy measure of $-X$, then
$\Pi^{(+)}$ for the restriction of $\Pi$ to $[0,\infty)$,
$\Pi^{(-)}$ for the restriction of  $\Pi^{\#}$ to $[0,\infty)$,
and
\begin{eqnarray}\label{tails}
\overline{\Pi}^{(+)}(x)&=&\Pi((x,\infty)),\nonumber\\
\overline{\Pi}^{(-)}(x)&=&\Pi
((-\infty,-x))\,,\nonumber\\
\overline{\Pi}(x)&=&\overline{\Pi}^{(+)}(x)+
\overline{\Pi}^{(-)}(x),\ x>0,
\end{eqnarray}
for the tails of $\Pi(\cdot)$.
Recall that we assume (\ref{nohold}) and that the L\'evy
measure $\Pi$ is  restricted to $[-1,1]$.
We will often make use of the L\'evy-It\^o decomposition, which
can be written as
\begin{equation}\label{eqlid}
X_t=\gamma t + \int_{[0,t]\times[-1,1]} x N({\rm d}s,{\rm
d}x)\,,
\qquad t\geq 0\,,
\end{equation}
where $N({\rm d}t,{\rm d}x)$ is a Poisson random measure on
$\R_+\times [-1,1]$
with intensity ${\rm d}t\Pi({\rm d}x)$ and the Poissonian
stochastic integral above is
taken in the compensated sense. See  Theorem 41 on page 31 in
\cite{prot} for details.

\section{Two-Sided Case}\label{s2}\
In this section we study 
two-sided crossings of power law boundaries at small times.
We wish to find a necessary and sufficient condition
for  \eqref{main}
for each value of $\kappa\ge 1/2$.
This question is completely answered in the next two theorems,
where the first can be viewed as a reinforcement of
(\ref{eqbg}).

\begin{thm}\label{th1}\ 
Assume \eqref{nohold}, and take $\kappa> 1/2$.  
When $X\in bv$, assume its drift is zero.

\noindent (i)\  If 
\begin{equation}\label{2}
\int_0^1 \overline{\Pi}(x^\kappa) {\rm d}x <\infty,
\end{equation}
then we have
\begin{equation}\label{1}
\lim_{t\downarrow 0}\frac{X_t}{t^\kappa}= 0 \ \hbox{ a.s.}
\end{equation}
(ii)\
Conversely, if \eqref{2} fails, then
\begin{equation*}
\limsup_{t\downarrow 0}\frac{|X_t-a(t)|}{t^\kappa}=\infty \hbox{
a.s.},
 \end{equation*}
for any deterministic function $a(t):[0,\infty)\mapsto \R$.
\end{thm}

 \begin{remark}\label{a} It is easy to check that
\eqref{2} is equivalent to
$$\int_{[-1,1]}|x|^{1/\kappa}\Pi({\rm d}x)<\infty.$$
The latter 
 holds for $0< \kappa\le 1/2$
for any L\'evy process, as a fundamental property of the L\'evy
canonical 
measure (\cite{b}, p.13).
\eqref{1} always holds  when $0<\kappa<1/2$,
as mentioned in Section \ref{s1},
but not necessarily when $\kappa=1/2$.
 \end{remark}

The case $\kappa=1/2$ which is excluded in Theorem \ref{th1}
turns out to have interesting and unexpected features. To put
these in context, 
let's first review some background.
Khintchine \cite{khin}
(see also Sato (\cite{satoa}, Prop. 47.12, p. 358)) 
showed that,
given any function $h(t)$, positive, continuous, and
nondecreasing 
in a neighbourhood of 0, and satisfying
$h(t) =o\left(\sqrt{t\log|\log t|}\right)$ as $t\downarrow 0$,
there is a L\'evy process with $\sigma^2=0$ such that
$\limsup_{t\downarrow 0}|X_t|/h(t)=\infty$ a.s.
For example, we can take
$h(t) =\sqrt{t} \left(\log|\log t|\right)^{1/4}$.
The corresponding L\'evy process satisfies
$\limsup_{t\downarrow 0}|X_t|/\sqrt{t}=\infty$ a.s.
Thus the implication \eqref{2} $\Rightarrow$ \eqref{1} 
 is not in general true when $\kappa=1/2$.

On the other hand, when $\kappa=1/2$,
Theorem \ref{th1} remains true for example when $X\in bv$, in
the
sense that both \eqref{2} and \eqref{1} then hold,
as follows from the fact that
$X_t=O(t)$ a.s. as $t\downarrow 0$.

Thus we can have
$\limsup_{t\downarrow 0}|X_t|/\sqrt{t}$ equal to $0$ or $\infty$
a.s., 
and we are led to ask  for a  
NASC
to decide between the alternatives.
We give such a condition in Theorem \ref{th2} and  furthermore
show that
$\limsup_{t\downarrow 0}|X_t|/\sqrt{t}$  may take a positive
finite value, a.s.
Remarkably, Theorem \ref{th2} simultaneously solves the
one-sided
problem.
These one sided cases are further investigated in Section
\ref{s3}, where it will be seen 
that, by contrast, the one and two sided situations are
completely different when $\kappa\ne 1/2$.

To state the theorem, we need the notation
\begin{equation}\label{Vdef}
V(x)= \int_{|y|\le x}y^2\Pi({\rm d}y), \ x>0.
\end{equation}

\begin{thm}\label{th2}
(The case $\kappa=1/2$.)
Assume \eqref{nohold}, and put 
$$
I(a )=\int_{0}^{1}x^{-1}\exp\left(-\frac{a ^2}{2V(x)}\right){\rm
d}x
$$
and
$$
\lambda_I^*:= \inf\{a >0: I(a )<\infty\}\in[0,\infty]
$$
(with the convention, throughout,
 that the inf of the empty set is $+\infty$).
Then, a.s.,
\begin{equation}\label{sqrt}
\displaystyle{
-\liminf_{t\downarrow 0}\frac{X_{t}}{\sqrt{t}}=
\limsup_{t\downarrow 0}\frac{X_{t}}{\sqrt{t}}=
\limsup_{t\downarrow 0}\frac{|X_{t}|}{\sqrt{t}}= \lambda_I^*.
}
\end{equation}
\end{thm}

\begin{remark}\label{r400}
(i)\ \eqref{sqrt} forms a nice counterpart to the iterated log
version in \eqref{khin} and \eqref{khin2}.

(ii)\ If \eqref{2} holds for some $\kappa>1/2$, 
then $V(x)=o(x^{2-1/\kappa})$ as
$x\downarrow 0$, so $\int_0^1 \exp\{-a^2/2V(x)\} {\rm d}x/x$
converges for all $\lambda>0$. Thus
$\lambda_I^*=0$ and
$\lim_{t\downarrow 0}X_t/\sqrt{t}=0$ a.s.  in this case,
according to Theorem \ref{th2}.
Of course, this agrees with Theorem \ref{th1}(i).

(iii)\
The convergence of
$\int_{|x|\le \e^{-e}} x^2\log|\log|x|| \Pi({\rm d}x)$ implies
the convergence of
$\int_0^1 \exp\{-a^2/2V(x)\} {\rm d}x/x$ for all $a>0$,
as is  easily checked,
hence we have
$\lim_{t\downarrow 0}|X_t|/\sqrt{t} =0$ a.s. for all
such L\'evy processes.
A finite positive value, a.s., for
$\limsup_{t\downarrow 0}|X_{t}|/\sqrt{t}$ can occur
only in  a small class
of L\'evy processes whose canonical measures have $\Pi({\rm
d}x)$ close 
to $|x|^{-3}{\rm d}x$ near $0$. For example,
we can find  a $\Pi$ such that, for small $x$,
$V(x)= 1/\log |\log x|$.
Then 
 $\int_0^{1/2} \exp\{-a^2/2V(x)\} {\rm d}x/x=
\int_0^{1/2} |\log x|^{-a^2/2} {\rm d}x/x$ 
is finite for $a>\sqrt{2}$ but infinite for $a\le \sqrt{2}$.
Thus $\limsup_{t\downarrow 0}|X_t|/\sqrt{t}=\sqrt{2}$ a.s.
for this process;
in fact, $\limsup_{t\downarrow 0}X_t/\sqrt{t}=\sqrt{2}$ a.s.,
and
 $\liminf_{t\downarrow 0}X_t/\sqrt{t}=-\sqrt{2}$ a.s.

(iv)\ 
Theorem \ref{th2} tells us that the only possible a.s. limit, 
as $t\downarrow 0$, of
$X_t/\sqrt{t}$ is 0, and that this occurs iff $\lambda_I^*=0$,
i.e., iff $I(\lambda)<\infty$ for all $\lambda>0$.
Similarly, the iterated log version in \eqref{khin2} gives that
the only possible a.s. limit, as $t\downarrow 0$, of 
$X_t/\sqrt{t\log|\log t|}$ is 0, and that this occurs iff $\sigma^2=0$.
When $\kappa>1/2$,
Theorem \ref{th1} gives that  $\lim_{t\downarrow 0}
X_t/t^\kappa=0$ a.s. iff $\int_0^1\pibar(x^\kappa){\rm
d}x<\infty$,
provided, when $\kappa\ge 1$, the drift $\delta=0$.

Another result in this vein is that we can have
$\lim_{t\downarrow 0}
X_t/t^\kappa=\delta$ a.s for a constant $\delta$ with
$0<|\delta|<\infty$, and $\kappa>0$, iff
$\kappa=1$, $X\in bv$, $\delta$ is the drift, and $\delta\ne 0$.

\end{remark}

The following corollary shows that centering has no effect in
the two-sided case.

\begin{cor}\label{cor1}
Assume \eqref{nohold}, and, if $X\in bv$, assume it has drift
zero.
Suppose
$\limsup_{t\downarrow 0}|X_t|/t^\kappa=\infty$ a.s., for some
$\kappa\ge 1/2$.
Then 
\begin{equation}\label{str}
\limsup_{t\downarrow 0}\frac{|X_t-a(t)|}{t^\kappa}=\infty\
{\rm  a.s.,\ for\  any\ nonstochastic}\ a(t).
\end{equation}
\end{cor}

Finally, in this section, 
Table 1 summarises the  conditions for \eqref{main}:
\bigskip


%
%
\message{S-Tables Macro v1.0, ACS, TAMU (RANHELP at VENUS.TAMU.EDU)}
%
%
\newhelp\stablestylehelp{You must choose a style between 0 and 3.}%
\newhelp\stablelinehelp{You should not use special hrules when stretching
a table.}%
\newhelp\stablesmultiplehelp{You have tried to place an S-Table inside another
S-Table.  I would recommend not going on.}%
%
%
\newdimen\stablesthinline
\stablesthinline=0.4pt
\newdimen\stablesthickline
\stablesthickline=1pt
%
%
\newif\ifstablesborderthin
\stablesborderthinfalse
\newif\ifstablesinternalthin
\stablesinternalthintrue
\newif\ifstablesomit
\newif\ifstablemode
\newif\ifstablesright
\stablesrightfalse
%
%
\newdimen\stablesbaselineskip
\newdimen\stableslineskip
\newdimen\stableslineskiplimit
%
%
\newcount\stablesmode
\newcount\stableslines
\newcount\stablestemp
\stablestemp=3
\newcount\stablescount
\stablescount=0
\newcount\stableslinet
\stableslinet=0
%
%
%
\newcount\stablestyle
\stablestyle=0
%
%
\def\stablesleft{\quad\hfil}%
\def\stablesright{\hfil\quad}%
%
%
\catcode`\@=\active%
%
%
\newcount\stablestrutsize
\newbox\stablestrutbox
\setbox\stablestrutbox=\hbox{\vrule height10pt depth5pt width0pt}
\def\stablestrut{\relax\ifmmode%
                         \copy\stablestrutbox%
                       \else%
                         \unhcopy\stablestrutbox%
                       \fi}%
%
%
\newdimen\stablesborderwidth
\newdimen\stablesinternalwidth
\newdimen\stablesdummy
\newcount\stablesdummyc
\newif\ifstablesin
\stablesinfalse
%
%
\def\begintable{\stablestart%
  \stablemodetrue%
  \stablesadj%
  \halign%
  \stablesdef}%
\def\begintableto#1{\stablestart%
  \stablemodefalse%
  \stablesadj%
  \halign to #1%
  \stablesdef}%
\def\begintablesp#1{\stablestart%
  \stablemodefalse%
  \stablesadj%
  \halign spread #1%
  \stablesdef}%
\def\stablesadj{%
  \ifcase\stablestyle%
    \hbox to \hsize\bgroup\hss\vbox\bgroup%
  \or%
    \hbox to \hsize\bgroup\vbox\bgroup%
  \or%
    \hbox to \hsize\bgroup\hss\vbox\bgroup%
  \or%
    \hbox\bgroup\vbox\bgroup%
  \else%
    \errhelp=\stablestylehelp%
    \errmessage{Invalid style selected, using default}%
    \hbox to \hsize\bgroup\hss\vbox\bgroup%
  \fi}%
\def\stablesend{\egroup%
  \ifcase\stablestyle%
    \hss\egroup%
  \or%
    \hss\egroup%
  \or%
    \egroup%
  \or%
    \egroup%
  \else%
    \hss\egroup%
  \fi}%
\def\stablestart{%
  \ifstablesin%
    \errhelp=\stablesmultiplehelp%
    \errmessage{An S-Table cannot be placed within an S-Table!}%
  \fi
  \global\stablesintrue%
  \global\advance\stablescount by 1%
  \message{<S-Tables Generating Table \number\stablescount}%
  \begingroup%
  \stablestrutsize=\ht\stablestrutbox%
  \advance\stablestrutsize by \dp\stablestrutbox%
  \ifstablesborderthin%
    \stablesborderwidth=\stablesthinline%
  \else%
    \stablesborderwidth=\stablesthickline%
  \fi%
  \ifstablesinternalthin%
    \stablesinternalwidth=\stablesthinline%
  \else%
    \stablesinternalwidth=\stablesthickline%
  \fi%
  \tabskip=0pt%
  \stablesbaselineskip=\baselineskip%
  \stableslineskip=\lineskip%
  \stableslineskiplimit=\lineskiplimit%
  \offinterlineskip%
  \def\borderrule{\vrule width \stablesborderwidth}%
  \def\internalrule{\vrule width \stablesinternalwidth}%
  \def\thinline{\noalign{\hrule height \stablesthinline}}%
  \def\thickline{\noalign{\hrule height \stablesthickline}}%
  \def\trule{\omit\leaders\hrule height \stablesthinline\hfill}%
  \def\ttrule{\omit\leaders\hrule height \stablesthickline\hfill}%
  \def\tttrule##1{\omit\leaders\hrule height ##1\hfill}%
  \def\stablesel{&\omit\global\stablesmode=0%
    \global\advance\stableslines by 1\borderrule\hfil\cr}%
  \def\el{\stablesel&}%
  \def\elt{\stablesel\thinline&}%
  \def\eltt{\stablesel\thickline&}%
  \def\elttt##1{\stablesel\noalign{\hrule height ##1}&}%
  \def\elspec{&\omit\hfil\borderrule\cr\omit\borderrule&%
              \ifstablemode%
              \else%
                \errhelp=\stablelinehelp%
                \errmessage{Special ruling will not display properly}%
              \fi}%
  \def\stmultispan##1{\mscount=##1 \loop\ifnum\mscount>3 \stspan\repeat}%
  \def\stspan{\span\omit \advance\mscount by -1}%
  \def\multicolumn##1{\omit\multiply\stablestemp by ##1%
     \stmultispan{\stablestemp}%
     \advance\stablesmode by ##1%
     \advance\stablesmode by -1%
     \stablestemp=3}%
  \def\multirow##1{\stablesdummyc=##1\parindent=0pt\setbox0\hbox\bgroup%
    \aftergroup\emultirow\let\temp=}
  \def\emultirow{\setbox1\vbox to\stablesdummyc\stablestrutsize%
    {\hsize\wd0\vfil\box0\vfil}%
    \ht1=\ht\stablestrutbox%
    \dp1=\dp\stablestrutbox%
    \box1}%
  \def\stpar##1{\vtop\bgroup\hsize ##1%
     \baselineskip=\stablesbaselineskip%
     \lineskip=\stableslineskip%
     \lineskiplimit=\stableslineskiplimit\bgroup\aftergroup\estpar\let\temp=}%
  \def\estpar{\vskip 6pt\egroup}%
  \def\stparrow##1##2{\stablesdummy=##2%
     \setbox0=\vtop to ##1\stablestrutsize\bgroup%
     \hsize\stablesdummy%
     \baselineskip=\stablesbaselineskip%
     \lineskip=\stableslineskip%
     \lineskiplimit=\stableslineskiplimit%
     \bgroup\vfil\aftergroup\estparrow%
     \let\temp=}%
  \def\estparrow{\vfil\egroup%
     \ht0=\ht\stablestrutbox%
     \dp0=\dp\stablestrutbox%
     \wd0=\stablesdummy%
     \box0}%
  \def@{\global\advance\stablesmode by 1&&&}%
  \def\@{\global\advance\stablesmode by 1&\omit\vrule width 0pt%
         \hfil&&}%
  \def\vt{\global\advance\stablesmode by 1&\omit\vrule width \stablesthinline%
          \hfil&&}%
  \def\vtt{\global\advance\stablesmode by 1&\omit\vrule width \stablesthickline%
          \hfil&&}%
  \def\vttt##1{\global\advance\stablesmode by 1&\omit\vrule width ##1%
          \hfil&&}%
  \def\vtr{\global\advance\stablesmode by 1&\omit\hfil\vrule width%
           \stablesthinline&&}%
  \def\vttr{\global\advance\stablesmode by 1&\omit\hfil\vrule width%
            \stablesthickline&&}%
  \def\vtttr##1{\global\advance\stablesmode by 1&\omit\hfil\vrule width ##1&&}%
  \stableslines=0%
  \stablesomitfalse}
\def\stablesdef{\bgroup\stablestrut\borderrule##\tabskip=0pt plus 1fil%
  &\stablesleft##\stablesright%
  &##\ifstablesright\hfill\fi\internalrule\ifstablesright\else\hfill\fi%
  \tabskip 0pt&&##\hfil\tabskip=0pt plus 1fil%
  &\stablesleft##\stablesright%
  &##\ifstablesright\hfill\fi\internalrule\ifstablesright\else\hfill\fi%
  \tabskip=0pt\cr%
  \ifstablesborderthin%
    \thinline%
  \else%
    \thickline%
  \fi&%
}%
\def\endtable{\advance\stableslines by 1\advance\stablesmode by 1%
   \message{- Rows: \number\stableslines, Columns:  \number\stablesmode>}%
   \stablesel%
   \ifstablesborderthin%
     \thinline%
   \else%
     \thickline%
   \fi%
   \egroup\stablesend%
\endgroup%
\global\stablesinfalse}
%
%


\def\stablesleft{\ \hfil}
\def\stablesright{\hfil \ }

\stablestyle=0

\centerline{\bf Table 1}

\bigskip

\begintable
Value of $\kappa$ @ NASC for $\limsup_{t\downarrow 0}
\vert X_t \vert/t^{\kappa} = \infty$ a.s. (when  $\sigma^2=0$)
\eltt
$0 \leq \kappa < \frac{1}{2}$ @ Never true 
\elt
$\kappa = \frac{1}{2}$ @ 
$\lambda_I^*=\infty$ (See Theorem \ref{th2})\elt
$\frac{1}{2} < \kappa \le 1$ @ $\int_0^1 \pibar(x^\kappa){\rm
d}x=\infty$ \elt
$\kappa > 1$, $X\in bv$, $\delta=0$ @ $\int_0^1
\pibar(x^\kappa){\rm d}x=\infty$ \elt
$\kappa > 1$, $X\in bv$, $\delta\ne 0$ @  Always true\elt
$\kappa > 1$, $X\notin bv$@  Always true
\endtable

\section{One-Sided Case}\label{s3}\
We wish to test for the finiteness or otherwise of
$\limsup_{t\downarrow 0}X_t/t^\kappa$,
so we proceed by finding conditions for
\begin{equation}\label{1s}
\limsup_{t\downarrow 0}\frac{X_{t}}{t^{\kappa}}=+\infty
\ \mathrm{a.s.}
\end{equation}

In view of the discussion in Section \ref{s1}, and the fact that
the case $\kappa=1/2$ is covered in Theorem \ref{th2},
we have only two cases to consider:

(a) $X\notin bv$, $1/2< \kappa<1$;

(b) $X\in bv$, with drift $\delta=0$, $\kappa>1$.

For Case (a), we need to define, for $1\geq y>0$, 
and  for $\lambda >0$,
$$
W(y):=\int_{0}^{y}\int_{x}^{1}z\Pi ^{(+)}({\rm d}z){\rm d}x, $$
and then
\begin{equation} \label{2.11}
J(\lambda):=\int_{0}^{1}\exp \left\{ -\lambda \left(
\frac{y^{{\frac{2\kappa
-1}{\kappa }}}}{W(y)}\right) ^{\frac{\kappa}{1-\kappa }}\right\}
\frac{{\rm d}y}{y}.
\end{equation}
Also let $\lambda_J^*:= \inf\{\lambda >0: J(\lambda )<\infty\}$.

\begin{thm}\label{thm3.1}
Assume \eqref{nohold} and 
keep $1/2< \kappa<1$. Then  \eqref{1s} holds if and only if
\newline
(i)\  $\int_0^1 \overline{\Pi}^{(+)}(x^\kappa){\rm d}x=\infty$,
or
\newline
(ii)\ 
$\int_0^1 \overline{\Pi}^{(+)}(x^\kappa){\rm d}x<\infty=\int_0^1
\overline{\Pi}^{(-)}(x^\kappa){\rm d}x$, and
$\lambda_J^*=\infty$.

When (i) and (ii) fail, we have in greater detail: suppose
\newline
(iii)\  $\int_0^1 \overline{\Pi}(x^\kappa){\rm d}x<\infty$,
or
\newline
(iv)\ $\int_0^1 \overline{\Pi}^{(+)}(x^\kappa){\rm d}x<\infty
=\int_0^1 \overline{\Pi}^{(-)}(x^\kappa){\rm d}x$
and $\lambda_J^*=0$.
 \newline
Then
\begin{equation}\label{ze}
\limsup_{t\downarrow 0}\frac{X_{t}}{t^{\kappa}}=0\ \mathrm{a.s.}
\end{equation}
Alternatively, suppose
\newline
(v)\
$\int_0^1 \overline{\Pi}^{(+)}(x^\kappa){\rm d}x<\infty
=\int_0^1 \overline{\Pi}^{(-)}(x^\kappa){\rm d}x$
and $\lambda_J^*\in(0,\infty)$.
Then 
\begin{equation}\label{ze2}
\limsup_{t\downarrow 0}\frac{X_{t}}{t^{\kappa}}=c \
\mathrm{a.s.,\ for\ some}\ c\in (0,\infty).
 \end{equation}
\end{thm}

\bigskip
\begin{remark}\label{r40}
(i)\ Using the integral criterion in terms of $J(\lambda)$, it's
easy to give examples of all three possibilities ($0$, $\infty$,
or in $(0,\infty)$) for
$\limsup_{t\downarrow 0}X_t/t^\kappa$, in the
situation of Theorem \ref{thm3.1}.

(ii)\
Note that $X\notin bv$ when 
$\int_0^1 \overline{\Pi}^{(+)}(x^\kappa){\rm d}x=\infty$
or
$\int_0^1 \overline{\Pi}^{(-)}(x^\kappa){\rm d}x=\infty$
in Theorem \ref{thm3.1},
because $ \overline{\Pi}^{(\pm)}(x^\kappa)\le
\overline{\Pi}^{(\pm)}(x)$
when $0<x<1$ and $\kappa<1$, so $\int_0^1 \overline{\Pi}(x){\rm
d}x=\infty$.

(iii)\ 
It may seem puzzling at first that
a second moment-like function, $V(\cdot)$, appears in Theorem
\ref{th2},
whereas $W(\cdot)$,  a kind of integrated  first moment
function,
appears in Theorem \ref{thm3.1}.
Though closely related, 
in general, $V(x)$ is not asymptotically equivalent to $W(x)$,
as $x\to 0$, and neither function is  asymptotically equivalent
to yet
another second moment-like function on $[0,\infty)$,
$U(x):= V(x)+x^2\pibar(x)$.
$V(x)$ arises naturally in the proof of  Theorem \ref{th2},
which uses a normal approximation to certain
probabilities,
whereas $W(x)$ arises naturally in the proof of  Theorem
\ref{thm3.1}, which uses Laplace transforms and works with
spectrally one-sided L\'evy processes.
It is possible to reconcile the different expressions; in fact,
Theorem \ref{th2} remains true if $V$ is replaced in the
integral
$I(\lambda)$ by $U$ or by $W$.
Thus these three functions are equivalent
in the context of Theorem \ref{th2} (but not in general).
We explain this in a little more detail following the proof of
Theorem 
\ref{thm3.1}.

\end{remark}

Next we turn to Case (b). When $X\in bv$ we can
define, for $0<x<1$,
 \begin{equation}\label{A-def}
A_+(x)= \int_0^x \overline{\Pi}^{(+)}(y){\rm d}y \
{\rm and}\ 
A_-(x)= \int_0^x \overline{\Pi}^{(-)}(y){\rm d}y.
 \end{equation}

\begin{thm}\label{thm3.2}
Assume \eqref{nohold},
suppose $\kappa>1$, $X\in bv$, and its drift $\delta=0$.
If
\begin{equation}\label{5.1}
 \int_{(0,1]}
\frac{\Pi^{(+)}({\rm d}x)}
{x^{-1/\kappa} +A_-(x)/x}=\infty
 \end{equation}
then \eqref{1s} holds.
Conversely, if \eqref{5.1} fails, then 
$\limsup_{t\downarrow 0}X_{t}/t^\kappa\le 0$ a.s.
\end{thm}

\bigskip
\begin{remark}\label{r4}
(i)\ It's natural to enquire whether 
\eqref{5.1} can be simplified by considering separately
integrals containing the components of the integrand in
\eqref{5.1}. This is not the case.
For each $\kappa>1$, it is possible to find a 
L\'evy process $X \in bv$
with drift $0$ for which \eqref{5.1} fails but
$$\int_{(0,1]} x^{\frac{1}{\kappa}}\Pi^{(+)}({\rm d}x)=\infty
= \int_{(0,1]} \left(x/A_-(x)\right)\Pi^{(+)}({\rm d}x).$$
The idea is to construct a continuous increasing concave
function which is linear on a sequence of intervals tending to
$0$, which can serve as an $A_-(x)$, and which oscillates around
the function $x\mapsto x^{1-1/{\kappa}}$.
Note that  \eqref{5.1} is equivalent to
$$\int_{(0,1]}  \min\left(x^{1/{\kappa}}, \frac{x}{A_-(x)}\right)
\Pi^{(+)}({\rm d}x)=\infty.$$ We will omit the details of the
construction.

(ii)\
It is possible to have $\limsup_{t\downarrow 0}X_{t}/t^\kappa<0$
a.s.,
in the situation of Theorem \ref{thm3.2},
when \eqref{5.1} fails; for example, when $X$ is the negative of
a subordinator with zero drift.
The value of the limsup can then be determined by applying Lemma
\ref{ronlem2} in Section \ref{s4}.

(ii)\
For another equivalence, we note that \eqref{5.1} holds
if and only if
 \begin{equation}\label{em}
\int_0^1
\overline{\Pi}^{(+)}(t^\kappa+X_t^{(-)})dt=\infty\ {\rm a.s.}
  \end{equation}
  where $X^{(-)}$ is a subordinator with drift $0$ and L\'evy
measure $\Pi^{(-)}$.
This can be deduced from
Erickson and Maller \cite{em}, Theorem 1, and provides a
connection
between the a.s. divergence of the L\'evy integral in
\eqref{em}
and the upcrossing condition \eqref{1s}.
\end{remark}

Table 2 summarises the  conditions for \eqref{1s}:

\bigskip
\centerline{\bf Table 2}
\bigskip
\begintable
Value of $\kappa$ @ NASC for 
$\limsup_{t\downarrow 0} X_t/t^{\kappa}=\infty$ a.s. 
(when  $\sigma^2=0$)
\eltt
$0 \leq \kappa < \frac{1}{2}$ @ Never true \elt
$\kappa = \frac{1}{2}$, $X\notin bv$ @ See Theorem \ref{th2}\elt
$\frac{1}{2} < \kappa < 1$, $X\notin bv$ @
See Theorem \ref{thm3.1} \elt
$\frac{1}{2} \le \kappa \le 1$, $X\in bv$@ Never true\elt
$\kappa > 1$, $X\in bv$, $\delta<0$ @  Never true\elt
$\kappa > 1$, $X\in bv$, $\delta=0$ @  See Theorem
\ref{thm3.2}\elt
$\kappa > 1$, $X\in bv$, $\delta>0$ @  Always true\elt
$\kappa \ge 1$, $X\notin bv$@  Always true
\endtable


\medskip

\vfill\eject
Our final theorem applies the foregoing results to give a
criterion for
\begin{equation}\label{2s}
\lim_{t\downarrow 0}\frac{X_{t}}{t^{\kappa}}=+\infty
\ \mathrm{a.s.}
\end{equation}
This is a stronger kind of divergence of the normed
process to $\infty$, for small times.
A straightforward analysis of cases, using our one and two-sided
results,  shows that \eqref{2s} never occurs if 
$0<\kappa\le 1$, if $\kappa>1$ and $X\notin bv$, or 
if $\kappa>1$ and  $X\in bv$  with negative drift.
If $\kappa>1$ and  $X\in bv$  with positive drift,
 \eqref{2s} always occurs.
That leaves just one case to consider, in:

\begin{thm}\label{thm3.3}
Assume \eqref{nohold},
suppose $\kappa>1$, $X\in bv$, and its drift $\delta=0$.
Then \eqref{2s} holds iff
\begin{equation}\label{33a}
K_X(d):=\int_{0}^{1}\frac{{\rm d}y}{y}\exp \left\{
- d
\frac{\left(A_+(y)\right)^{\frac{\kappa}{\kappa-1}}}{y}\right\}
<\infty,\ {\rm for \ all}\ d>0,
 \end{equation}
and
\begin{equation}\label{33b}
\int_{(0,1]} \frac{x } {A_+(x)}\Pi^{(-)}({\rm d}x) <\infty.
\end{equation}
\end{thm}

\noindent {\bf Concluding Remarks.}\ 
It's natural to enquire about a one-sided version of Corollary
\ref{cor1}: when is
\begin{equation}\label{str2}
\limsup_{t\downarrow 0}\frac{X_t-a(t)}{t^\kappa}<\infty\
{\rm  a.s.,\ for\  some\ nonstochastic}\ a(t)?
\end{equation}
Phrased in such a general way the question is not interesting
since we can always make $X_t=o(a(t))$ a.s as $t\downarrow 0$ by
choosing $a(t)$ large enough by comparison with $X_t$
(e.g., $a(t)$ such that $a(t)/\sqrt{t\log|\log t|}\to \infty$,
as $t\downarrow 0$, will do, by \eqref{khin}),
so the limsup in \eqref{str2} becomes negative.
So we would need to restrict $a(t)$
in some way. Section \ref{s3} deals with the case $a(t)=0$.
Another choice is to take $a(t)$ as a natural centering function
such as $EX_t$ or as a median of $X_t$.
However, in our small time situation, $EX_t$ is  essentially 
0 or the drift of $X$, so we are led back to the case $a(t)=0$
again (and similarly for the median).
Of course there may be other interesting choices of $a(t)$ in
some applications, and there is the wider issue 
of replacing $t^\kappa$ by a more general norming function. Some
of our results in Sections \ref{s4} and \ref{s5} address the
latter, but we will not pursue these points further here.

\section{Proofs for Section \ref{s2}}\label{s4} 
\subsection{Proof of Theorem  \ref{th1}}
The proof relies on a pair of technical
results which we will establish first.
Recall the notation  $V(x)$ in \eqref{Vdef}.

\begin{prop}\label{prop2}
Let $b: \R_+\to [0,\infty)$ be any
non-decreasing function such that
$$\int_0^1 \overline\Pi(b(x)){\rm d}x<\infty
\quad\hbox{and}\quad
\int_0^1V(b(x)) b^{-2}(x) {\rm d}x <\infty.$$
Then
$$\limsup_{t\downarrow 0}\frac{|X_t-a(t)|}{b(4t)}\leq 1\ {\rm a.s.,}$$
where
\begin{equation}\label{adef}
a(t):=\gamma t - \int_0^t {\rm d}s\int_{b(s)<
|x|\leq 1}x\Pi({\rm d}x),\ t\geq 0.
\end{equation}
\end{prop}

\bigskip \noindent{\it Proof of Proposition \ref{prop2}:}\
Recall the L\'evy-It\^o decomposition (\ref{eqlid}).
In this setting, it is convenient to introduce 
$$X^{(1)}_t:= \int_{[0,t]\times[0,1]}
{\bf 1}_{\{|x|\leq b(s)\}} x N({\rm d}s,{\rm d}x)$$
and
$$X^{(2)}_t:=\gamma t + \int_{[0,t]\times[0,1]}
{\bf 1}_{\{b(s)< |x|\leq 1\}} x N({\rm d}s,{\rm d}x)\,,$$
where again the stochastic integrals are taken in the
compensated sense.
Plainly, $X=X^{(1)}+X^{(2)}$.

The assumption $\int_0^1 \overline\Pi(b(x)){\rm d}x<\infty$
implies that
 $$N\left(\{(s,x): 0\leq s \leq t \hbox{ and } b(s)< |x|\leq
1\}\right)=0$$
 whenever $t>0$ is sufficiently small a.s., and in this
situation $X^{(2)}$ is just $\gamma t$ minus the compensator,
a.s.;
i.e., 
$$X^{(2)}_t=\gamma t - \int_0^t {\rm d}s\int_{b(s)< |x|\leq
1}x\Pi({\rm d}x)=a(t).$$

On the one hand, $X^{(1)}$ is a square-integrable martingale
with oblique bracket 
$$\langle X^{(1)}\rangle_t=
\int_0^t {\rm d} s \int_{|x|\leq b(s)}x^2\Pi({\rm d}x)
=\int_0^t V(b(s)){\rm d} s \leq tV(b(t)).$$
By Doob's maximal inequality, we have for every $t\geq 0$
$$P(\sup_{0\leq s\leq t}|X^{(1)}_s|> b(2t))
\leq 4 t V(b(t)) b^{-2}(2t).$$
On the other hand, the assumptions that $b(t)$ is non-decreasing
and that
$\int_0^1{\rm d}xV(b(x))b^{-2}(x)<\infty$
 entail
$$\sum_{n=1}^{\infty} 2^{-n}V(b(2^{-n}))
b^{-2}(2^{-n+1})<\infty.$$
By the Borel-Cantelli lemma, we thus see that 
$$\lim_{n\to \infty}\frac{\sup_{0\leq s\leq
2^{-n}}|X^{(1)}_s|}{b(2^{-n+1})}\leq 1
\qquad \hbox{a.s.},$$
and the proof is completed by a standard argument of
monotonicity.
\hfill\bbox

\begin{prop}\label{prop1}
Suppose there are deterministic  functions $a:\R_+\to \R$ and
$b:
(0,\infty)\to (0,\infty)$, with $b$ measurable,
such that
\begin{equation}\label{3a}
P\left( \limsup_{t\downarrow 0}
{\frac{|X_t-a(t)|}{b(t)}} < \infty\right) >0.
\end{equation}
Then there is some finite constant $C$ such that
\begin{equation}\label{3}
\int_0^1 \overline{\Pi}(Cb(x)) {\rm d}x <\infty.
\end{equation}
\end{prop}

\bigskip \noindent{\it Proof of Proposition \ref{prop1}:}
Symmetrise $X$ by subtracting an independent equally
distributed $X'$ to get $X_t^{(s)}=X_t-X_t'$, $t\ge 0$. Then
\eqref{3a}  and Blumenthal's $0$-$1$ law imply
there is some finite constant $C$ such that
\begin{equation}\label{eq2}
\limsup_{t\downarrow 0}\frac{|X^{(s)}_{t}|}{b(t)}<\frac{C}{2},\
\hbox{a.s.}
\end{equation}

Suppose now that (\ref{3}) fails. Note that 
$\overline{\Pi}^{(s)}(\cdot)=2\overline{\Pi}(\cdot)$,
where $\Pi^{(s)}$ is the L\'evy measure of $X^{(s)}$, so that
 $\int_0^1 \overline{\Pi}^{(s)}(Cb(x)) {\rm d}x =\infty$.
Then from the L\'evy-It\^{o}
decomposition, we have that,
for every $\varepsilon >0$, 
$$
\#\{t\in (0,\varepsilon]: |\Delta_t^{(s)}|>Cb(t)\}=\infty,\
{\rm a.s.,}
$$
where $\Delta_t^{(s)}=X_t^{(s)}-X_{t-}^{(s)}$.
But whenever $|\Delta_t^{(s)}|>Cb(t)$, we must have 
$|X_{t-}^{(s)}|>Cb(t)/2$ or $|X_{t}^{(s)}|>Cb(t)/2$; which
contradicts
(\ref{eq2}). Thus \eqref{3} holds.
\hfill\bbox

\medskip 

Finally, we will need an easy deterministic bound.

\begin{lem}\label{Ldet} Fix some $\kappa\ge 1/2$ and assume
\begin{equation}\label{2-bis}
\int_{|x|<1}|x|^{1/\kappa}\Pi({\rm d}x)<\infty.
\end{equation}
When $\kappa\ge 1$, $X\in bv$ and we suppose further that the
drift coefficient $\delta=\gamma -\int_{|x|\le 1} x\Pi({\rm
d}x)$ is $0$. Then, as $t\to 0$,
$$a(t)=\gamma t - \int_0^t {\rm d}s\int_{s^{\kappa}<|x|\leq 
1}x\Pi({\rm d}x)=o(t^{\kappa}).$$
\end{lem}

\bigskip \noindent{\it Proof of Lemma \ref{Ldet}:}\
Suppose first
$\kappa<1$. For every $0<\varepsilon<\eta<1$, we have
\begin{eqnarray*}
\int_{\varepsilon< |x| \leq 1}|x|\Pi({\rm d}x) &\leq&
\varepsilon^{1-1/\kappa}
\int_{ |x|\leq \eta}|x|^{1/\kappa}\Pi({\rm d}x)+ \int_{\eta<
|x|\leq 1}|x|\Pi({\rm d}x)\\
&=&  \varepsilon^{1-1/\kappa}o_\eta+c(\eta), \ {\rm say,}
\end{eqnarray*}
where, by (\ref{2-bis}), $\lim_{\eta\downarrow 0}o_\eta=0$.
Since $\kappa<1$, it follows that 
$$\limsup_{t\downarrow 0}|a(t)|t^{-\kappa}\leq \kappa^{-1}o_\eta,$$
and as we can take $\eta$ arbitrarily small,
 we conclude that $a(t)=o(t^{\kappa})$.

In the case $\kappa \geq 1$, $X$ has bounded variation with zero
drift coefficient. We may rewrite
$a(t)$ in the form
$$a(t)=\int_0^t {\rm d}s\int_{ |x|\leq s^{\kappa}}x\Pi({\rm
d}x).$$
The assumption (\ref{2-bis}) entails
$\int_{ |x|\leq \varepsilon}|x|\Pi({\rm
d}x)=o(\varepsilon^{1-1/\kappa})$
and we again conclude that $a(t)=o(t^{\kappa})$.
\hfill \bbox

\medskip

We now have all the ingredients to establish Theorem \ref{th1}.

\bigskip \noindent{\it Proof of Theorem \ref{th1}:}\ 
Keep $\kappa>1/2$ throughout. (i) Suppose
(\ref{2}) holds, which is equivalent to (\ref{2-bis}).
 Writing $|x|^{1/\kappa}=|x|^{1/\kappa-2}x^2$,
we see from an integration by parts (Fubini's theorem) that
$\int_0^1 V(x) x^{1/\kappa-3}{\rm d}x<\infty$. Note that the
assumption that $\kappa\neq 1/2$ is crucial in this step. The
change of variables $x=y^{\kappa}$
now gives that $\int_0^1V(y^{\kappa})y^{-2\kappa}{\rm d}
y<\infty$.
We may thus apply Proposition \ref{prop2} and get that
$$\limsup_{t\downarrow 0}\frac{|X_t-a(t)|}{t^{\kappa}}\leq 4^{\kappa}
\qquad \hbox{a.s.}$$
where $a(t)$ is as in Lemma \ref{Ldet}.
We thus have shown that when
(\ref{2}) holds, 
$$\limsup_{t\downarrow 0}\frac{|X_t|}{t^{\kappa}}\leq 4^{\kappa}
\qquad \hbox{a.s.}$$
For every $c>0$, the time-changed process $X_{ct}$ is a L\'evy
process
with L\' evy measure $c\Pi$, so we also have
$\limsup_{t\downarrow 0}{|X_{ct}|}t^{-\kappa}\leq 4^{\kappa}$ a.s.
As we may take $c$ as large as we wish, we conclude that
$$\lim_{t\downarrow 0}\frac{|X_t|}{t^{\kappa}}=0 \qquad \hbox{a.s.}$$
 
(ii) By Proposition \ref{prop1}, if
$$ P\left( \limsup_{t\downarrow 0}
{\frac{|X_t-a(t)|}{b(t)}} < \infty\right) >0,$$
then $\int_0^1\overline{\Pi}(Cx^{\kappa}){\rm d} x<\infty$
for some finite constant $C$. By an obvious change of variables,
this shows that (\ref{2}) must hold. 
This completes the proof of Theorem \ref{th1}.  \hfill\bbox

\medskip
Finally we establish Corollary \ref{cor1}.

\bigskip \noindent{\it Proof of Corollary \ref{cor1}:}\
This hinges on the fact that when \eqref{str} fails, then
$(X_t-a(t))/t^\kappa \topr 0$,
with $a(t)$ defined as in \eqref{adef} -- even for $\kappa=1/2$.
We will omit the details.
\hfill\bbox

\subsection{Proof of Theorem \ref{th2}}
We now turn our attention to Theorem \ref{th2} and develop some
notation and material in this direction.
Write, for $b>0$, 
\begin{equation}\label{2.2z}
X_t=Y^{(b)}_t+Z^{(b)}_t,
\end{equation}
 with
\begin{eqnarray}
Y^{(b)}_t&:=&\int_{[0,t]\times[-1,1]}{\bf 1}_{\{|x|\leq
b\}}xN({\rm d}s,{\rm d}x)\,,
\label{2.2a}\\
Z^{(b)}_t&:=&\gamma t+\int_{[0,t]\times[-1,1]}{\bf
1}_{\{b<|x|\}}x N({\rm d}s,{\rm d}x)\,,
\nonumber
\end{eqnarray}
where 
$N({\rm d}s,{\rm d}x)$ is a Poisson random measure on
$[0,\infty)\times [-1,1]$
with intensity ${\rm d}s \Pi({\rm d}x)$, and the stochastic
integrals are taken in
the compensated sense.

\begin{lem}\label{LR0}\ (No assumptions on $X$.)
For every $0<r<1$ and $\varepsilon >0$,
we have
$$\sum_{n=1}^{\infty}P(\sup_{0\leq t \leq
r^n}|Z^{(r^{n/2})}_{t}| >\varepsilon r^{n/2})<\infty,$$
and as a consequence,
$$\lim_{n\to\infty} r^{-n/2}\sup_{0\leq t \leq
r^n}|Z^{(r^{n/2})}_{t}| =0\ {\rm a.s.}$$
\end{lem}

\bigskip \noindent{\it Proof of Lemma \ref{LR0}:}\
Introduce, for every integer $n$, the set
$$A_n:=[0,r^n]\times([-1,-r^{n/2})\cup(r^{n/2},1])\,,$$
so that 
\begin{eqnarray*}
\sum_{n=1}^{\infty}P(N(A_n)>0)&\leq &\sum_{n=1}^{\infty}
r^n\overline{\Pi}(r^{n/2})\\
&\leq&  (1-r)^{-1}\sum_{n=1}^{\infty}\int^{r^n}_{r^{n+1}}
\pibar(\sqrt{x}){\rm d}x\\
&\leq&  (1-r)^{-1}\int_{0}^{1}
\pibar(\sqrt{x}){\rm d}x.
\end{eqnarray*}
As the last integral is finite (always), we have from the
Borel-Cantelli lemma
that $N(A_n)=0$ whenever $n$ is sufficiently large, a.s.

On the other hand, on the event $N(A_n)=0$, we have
$$Z^{(r^{n/2})}_{t}=t\left(\gamma -\int_{r^{n/2}<|x|\leq
1}x\Pi({\rm d}x)\right),\qquad 0\leq t \leq r^n.$$
Again as a result of the convergence of $\int_{|x|\le 1} x^2 \Pi({\rm d}x)$,
the argument in Lemma \ref{Ldet} shows that the supremum over
$0\leq t \leq r^n$ 
of the absolute value of the right-hand side
is $o(r^{n/2})$. The Borel-Cantelli lemma completes the proof.
\hfill\bbox

In view of Lemma \ref{LR0} we can concentrate on 
$Y_t^{(\sqrt t)}$ in (\ref{2.2z}).
We next prove:

\begin{lem}\label{ron2lem1}
Let $Y$ be a L\'evy process with canonical measure
$\Pi_Y$, satisfying $EY_1=0$ and $m_4<\infty$,
where $m_k:= \int_{x\in \R}|x|^k\Pi_Y ({\rm d}x)$,
$k=2,3,\cdots$.
\newline (i)\ Then 
\[
\lim_{t\downarrow 0} \frac{1}{t}E|Y_t|^3=m_3.
\]
\newline (ii)\
For any $x>0$, $t>0$,  we have the bound
\begin{equation}\label{be1}
\left|P(Y_t>x\sqrt{tm_2})-\overline{F}(x)\right|
\leq \frac{Am_3}{\sqrt{t} m_2^{3/2}(1+x)^3},
\end{equation}
where 
$$\overline{F}(x)=\int_x^\infty \e^{-y^2/2}{\rm
d}y/\sqrt{2\pi}=\frac{1}{2}{\rm erfc}(x/\sqrt2)$$
is the tail of the standard normal
distribution function, and $A$ is an absolute constant.
\end{lem}

\bigskip \noindent{\it Proof of Lemma \ref{ron2lem1}:}\
(i)\ We can calculate
$EY_t^4=tm_4+3t^2m_2^2$. So by Chebychev's inequality for second
and fourth moments, for $x>0$, $t>0$,
\[
\frac{1}{t}P(|Y_t|>x) \le \frac{m_2}{x^2}
 \boldsymbol{1}_{\{0<x\le 1\}}
 +\frac{m_4+3tm_2^2}{x^4} \boldsymbol{1}_{\{x> 1\}}.
\]
We can also calculate
\[
\frac{1}{t}E|Y_t|^3=\frac{3}{t} \int_0^\infty x^2P(|Y_t|>x) {\rm
d}x.
\]
By \cite{b}, Ex. 1, p. 39, 
$P(|Y_t|>x)/t\to \pibar_Y(x)$, as $t\downarrow 0$, for each
$x>0$.
The result (i) follows by dominated convergence.

(ii)\ 
 Write 
$Y_t = \sum_{i=1}^n Y(i,t)$, for $n=1,2,\cdots$, where
$Y(i,t):=Y(it/n)-Y((i-1)t/n)$ are i.i.d., each with the
distribution of $Y(t/n)$. 
According to  a non-uniform Berry-Esseen bound
(Theorem 14, p.125 of Petrov \cite{pet}),  for each
$n=1,2,\cdots,$ \eqref{be1} holds with the righthand side
replaced by
\[
\frac{AE|Y(t/n)|^3}{\sqrt{n}(tm_2/n)^{3/2}(1+x)^3}
=
\frac{AE|Y(t/n)|^3/(t/n)}{\sqrt{t}m_2^{3/2}(1+x)^3}.
\]
By Part (i) this tends as $n\to \infty$ to  the righthand side
of  \eqref{be1}.  \hfill\bbox

\begin{prop}\label{ron2prop1}
In the notation \eqref{2.2a},
we have, for $a>0$, $0<r<1$,
\begin{equation}
\sum_{n\ge 0}
 P\left(Y_{r^n}^{(r^{n/2})}>ar^{n/2}\right)<\infty 
\Longleftrightarrow 
\int_0^1
 \sqrt{V(x)}
\exp\left( \frac{-a^2}{2V(x)}\right)\frac{{\rm d}x}{x}
<\infty.
\label{be3}
\end{equation}
\end{prop}

\bigskip \noindent{\it Proof of Proposition \ref{ron2prop1}:}\
For every fixed $t>0$, 
$Y_s^{(\sqrt{t})}$ is the compensated sum of jumps of $X$
smaller in magnitude than $\sqrt{t}$, up to time $s$.
It is a centered L\'evy process with 
canonical measure
$\boldsymbol{1}_{\{|x|\le \sqrt{t}\}}\Pi({\rm d}x)$, $x\in \R$.
Applying  Lemma \ref{ron2lem1}, we get
$m_2=V(\sqrt{t})$
and
$m_3=\int_{|y|\le  \sqrt{t}}|y|^3\Pi ({\rm d}y)=\rho(
\sqrt{t})$,
say.
Then we get, for $x>0$,
\begin{equation*}
|P\left(Y_t^{(\sqrt{t})}>x\sqrt{tV(\sqrt{t})}\right)-\overline{F}(x)|\leq 
\frac{A\rho (\sqrt{t})}{\sqrt{tV^3(\sqrt{t})}(1+x)^3}.
\end{equation*}
Replacing $x$ by
$a/\sqrt{V(\sqrt{t})}$, $a>0$,
we have 
\[
\left|P\left(Y_t^{(\sqrt{t})}>a\sqrt{t}\right)
-\overline{F}\left(a/\sqrt{V(\sqrt{t})}\right)\right| \leq
\varepsilon (t):=
\frac{A\rho(\sqrt{t})}{\sqrt{t}a^{3}},
\]
and we claim that $\sum \varepsilon (r^{n})<\infty.$ In fact,
for some $c>0$,
\begin{eqnarray*}
\sum_{n=0}^\infty\frac{\rho(r^{n/2})}{r^{n/2}} 
&=&\sum_{n=0}^\infty 
\frac{1}{r^{n/2}}
\sum_{j\geq n}
\int_{r^{(j+1)/2}<|y|\leq r^{j/2}}|y|^3\Pi ({\rm d}y) \\
&=&
\sum_{j=0}^\infty \left( \sum_{n=0}^{j} r^{-n/2}\right) 
\int_{r^{(j+1)/2}<|y|\leq r^{j/2}}|y|^{3}\Pi
({\rm d}y) \\
&\leq &
c\sum_{j=0}^\infty r^{-j/2} \int_{r^{(j+1)/2}<|y|\leq
r^{j/2}}|y|^{3}\Pi ({\rm d}y)%
 \\
&\leq &
c\sum_{j=0}^\infty
\int_{r^{(j+1)/2}<|y|\leq r^{j/2}}y^{2}\Pi ({\rm d}y)
\\
&=&c\int_{|y|\leq 1}y^{2}\Pi ({\rm d}y)<\infty.
\end{eqnarray*}
The result (\ref{be3}) follows, since the monotonicity of
$\overline{F}$
shows that the convergence of
$\sum_{n\ge 1} \overline{F}\left(a/\sqrt{V(r^{n/2})}\right)$
 is equivalent to that of
\begin{equation*}
\int_{0}^{1}\overline{F}\left(a/\sqrt{V(\sqrt{x})}\right)\frac{{\rm
d}x}{x}
=\sum \int_{r^{
\frac{n+1}{2}}}^{r^{\frac{n}{2}}}\overline{F}\left(a/\sqrt{V(\sqrt{x})}\right)
\frac{{\rm d}x}{x},
\end{equation*}
and it is well-known that 
$\overline{F}(x)\backsim (2\pi)^{-1/2}x^{-1}\e^{-x^{2}/2}$ as
$x\to \infty$.
\hfill\bbox

We can now establish Theorem \ref{th2}.

\bigskip \noindent{\it Proof of Theorem \ref{th2}:}\
Recall the definition of  $I(\cdot)$ in the statement of 
Theorem \ref{th2}. We will first show
that for every given  $ a>0$
 \begin{equation}\label{up}
I(a)<\infty\ \ \Rightarrow\
\limsup_{t\downarrow 0}\frac{X_t}{\sqrt{t}}\le a\ {\rm a.s.}
\end{equation}
To see this, observe when $I(a)<\infty$, the integral in
\eqref{be3} converges, hence so
does the series.
Use the maximal inequality in 
Theorem 12, p.50 of Petrov \cite{pet},
to get, for $t>0$, $b>0$, $x>0$,
\begin{align}\label{5st}
P\left(\sup_{0<s\le t}Y_s^{(b)}>x\right)
&=\lim_{k \rightarrow \infty} 
P\left(\max_{1\le j\le \lceil kt \rceil}
Y_{j/k}^{(b)}>x\right)
 \nonumber \\
&=\lim_{k \rightarrow \infty} 
P\left(\max_{1\le j\le \lceil kt \rceil}
\sum_{i=1}^j \Delta(i,k,b)>x\right)
 \nonumber \\
&\le \limsup_{k \rightarrow \infty}
2 P\left(\sum_{i=1}^{\lceil kt \rceil} 
\Delta(i,k,b)>x- \sqrt{2ktV(b)/k} \right)
 \nonumber \\
&= 2 P\left(Y_t^{(b)}>x-\sqrt{2tV(b)}
\right),
\end{align}
where we note that
$\left(\Delta(i,k,b):=Y_{i/k}^{(b)}- Y_{(i-1)/k}^{(b)}\right)_{i\ge 1}$
are i.i.d., each with expectation 0 and variance equal to
$V(b)/k$.
Given $\veps>0$, replace $t$ by $r^n$, $b$ by $r^{n/2}$, 
and $x$ by  $ar^{n/2}+\sqrt{2r^nV(r^{n/2})}$,
which is not larger than
$(a+\veps)r^{n/2}$, once $n$ is large enough, in \eqref{5st}.
The convergence of the series in \eqref{be3} then gives
\[
\sum_{n\ge 0}
P\left(\sup_{0<s\le
r^n}Y_s^{(r^{n/2})}>(a+\veps)r^{n/2}\right)<\infty\ {\rm for\ all}\ \veps>0.
\]
Hence by the Borel-Cantelli Lemma
\begin{equation}\label{6st}
\limsup_{n \rightarrow \infty}
\frac{\sup_{0<t\le r^n}Y_t^{(r^{n/2})}}{r^{n/2}}\le a,\ 
{\rm a.s.}
\end{equation}
Using \eqref{2.2z},
together with Lemma \ref{LR0} and \eqref{6st}, gives
\[
\limsup_{n \rightarrow \infty}
\frac{\sup_{0<t\le r^n}X_t}{r^{n/2}}\le a,\ {\rm a.s.}
\]
By an argument of monotonicity, this yields
\[
\limsup_{t \downarrow 0}
\frac{X_t}{\sqrt{t}}
\le \frac{a}{\sqrt{r}},\ {\rm a.s.}
\]
Then  let $r \uparrow 1$ to get
$\limsup_{t\downarrow 0}X_t/\sqrt{t}\le a$ a.s.

For a reverse inequality, we show that for every $a>0$,
 \begin{equation}\label{down}
I(a)=\infty\ \Rightarrow \
\limsup_{t\downarrow 0}\frac{X_t}{\sqrt{t}}\ge a,\ {\rm a.s.}
\end{equation}
To see this, suppose that $I(a)=\infty$ for a given  $ a>0 $.
Then the integral in \eqref{be3} diverges when $a$ is replaced
by
$a-\veps$ for an arbitrarily small
$\veps>0$, because
$V(x)\ge \veps\exp(-\veps/2V(x))/2$, for $\veps>0$, $x>0$.
Hence, keeping in mind \eqref{2.2z}, Lemma \ref{LR0} and
Proposition
 \ref{ron2prop1}, we deduce
\begin{equation}\label{rev}
\sum_{n\ge 0} P\left(X_{r^n}>a'r^{n/2}\right)=\infty
\end{equation}
for all $a'<a$.
For a given $\veps>0$, define for every integer $n\geq 0$ the
events
\begin{eqnarray*}
A_n&=&\{X_{r^n/(1-r)}-X_{r^{n+1}/(1-r)}>a' r^{n/2}\},\\
B_n&=&\{|X_{r^{n+1}/(1-r)}|\le \veps r^{n/2}\}.
\end{eqnarray*}
Then the $\{A_n\}_{n\ge 0}$ are independent, and each $B_n$ is
independent of
the collection
$\{A_n, A_{n-1}, \cdots, A_0\}$.
Further, $\sum_{n\ge 0}P(A_n)=\infty$ by \eqref{rev}, so
$P(A_n\ {\rm i.o.})=1$.
It can be deduced easily from \cite{b}, Prop. 2(i), p.16,
that $X_t/\sqrt{t} \topr 0$, as
$t\downarrow 0$, since \eqref{nohold} is enforced.
Thus $P(B_n)\to 1$ as $n\to\infty$, and then, by the
Feller--Chung lemma (\cite{CT}, p. 69)
we can deduce that $P(A_n\cap B_n \ {\rm i.o.})=1$.
This implies
$P(X_{r^n/(1-r)}> (a'-\veps)r^{n/2}\ {\rm i.o.})=1$,
thus
\[
\limsup_{t\downarrow 0}
\frac{X_t}{\sqrt{t}} \ge (a'-\veps)\sqrt{1-r}\ {\rm a.s.},
\]
in which we can let $a'\uparrow a$,
$\veps\downarrow 0$ and $r\downarrow 0$ to get \eqref{down}.

As just mentioned, we have $X_t/\sqrt{t} \topr 0$, as
$t\downarrow 0$, so
$\liminf_{t\downarrow 0} X_t/\sqrt{t}  \le 0 
\le \limsup_{t\downarrow 0} X_t/\sqrt{t} $ a.s.
Together with \eqref{up} and \eqref{down}, this gives the
statements in
Theorem \ref{th2}
(replace $X$ by $-X$ to deduce the liminf statements from the
limsup,
noting that this leaves $V(\cdot)$ unchanged).  \hfill\bbox
 
\section{Proofs for Section \ref{s3}}\label{s5}

\subsection{Proof of Theorem \ref{thm3.1}}\
We start with some notation and technical results.
Recall we assume (\ref{nohold}).

 Take $0<\kappa<1$  and suppose first
\begin{equation}\label{4.1}
 \int_0^1 \overline{\Pi}^{(+)}(x^\kappa){\rm d}x=\infty.
\end{equation}
Define, for $0<x<1$,
\[
\rho_\kappa(x)= \frac{1}{\kappa}
\int_x^1 y^{\frac{1}{\kappa}-1} \overline{\Pi}^{(+)}(y){\rm d}y
=
\int_{x^\frac{1}{\kappa}}^1\overline{\Pi}^{(+)}(y^\kappa){\rm
d}y.
\]
Since $\overline{\Pi}^{(+)}(x)>0$ for all small $x$,
$\rho_\kappa(x)$ is strictly decreasing in a neighbourhood of 0,
thus 
$x^{-1/\kappa} \rho_\kappa(x)$
is also  strictly decreasing in a neighbourhood of 0,
and tends to $\infty$ as $x\downarrow 0$
(because  of \eqref{4.1}).
Also define
\[
U_+(x)= 2 \int_0^xy \overline{\Pi}^{(+)}(y){\rm d}y.
\]
By a similar argument, $x^{-2}U_+(x)$
is strictly decreasing in a neighbourhood of 0,
and tends to $\infty$ as $x\downarrow
0$.

Next, given $\alpha\in (0,\kappa)$,  define, for $t>0$,
\[
c(t)=\inf\left\{x>0: 
\rho_\kappa(x)x^{-1/\kappa}+ x^{-2}U_+(x)\le \alpha/t\right\}.
\]
Then $0<c(t)<\infty$ for $t>0$, $c(t)$ is strictly increasing,
$\lim_{t\downarrow 0}c(t)=0$,  and
\begin{equation}\label{4.7}
\frac{t\rho_\kappa(c(t))}{c^{\frac{1}{\kappa}}(t)}
+\frac{tU_+(c(t))}{c^2(t)}= \alpha.
\end{equation}
Since
 $\lim_{t\downarrow 0} \rho_\kappa(t)=\infty$, we have
 $\lim_{t\downarrow 0}c(t)/t^\kappa=\infty$.

We now point out that \eqref{4.1} can be reinforced as follows.

\begin{lem} \label{L4.3}
 The condition \eqref{4.1} implies that
 $$\int_0^1\pibar^{(+)}(c(t)){\rm d}t=\infty\,.$$
 \end{lem}
 
 \bigskip \noindent{\it Proof of Lemma \ref{L4.3}:}\
 We will first establish
 \begin{equation}\label{4.3}
\int_0^1 \frac{\Pi({\rm d}x)}{x^{-1/\kappa} \rho_\kappa(x)+
x^{-2}U_+(x)}
=\infty.
 \end{equation}
Suppose \eqref{4.3} fails.
Since 
\[
f(x):=\frac{1}{x^{-1/\kappa} \rho_\kappa(x)+ x^{-2}U_+(x)}
\]
is nondecreasing (in fact,  strictly increasing in a
neighbourhood of 0) with $f(0)=0$, 
for every $\veps>0$ there is an $\eta>0$ such that, for all
$0<x<\eta$,
\[
\veps\ge \int_x^\eta f(z)\Pi({\rm d}z)
\ge
f(x)\left(\overline{\Pi}^{(+)}(x)-\overline{\Pi}^{(+)}(\eta)\right),
\]
giving
\[
f(x)\overline{\Pi}^{(+)}(x)\le \veps+
f(x)\overline{\Pi}^{(+)}(\eta)
=\veps+o(1), \ {\rm as}\ x \downarrow 0.
\]
Letting ${\veps\downarrow 0}$ shows that
\begin{equation}\label{4.4}
\lim_{x\downarrow 0} f(x)\overline{\Pi}^{(+)}(x)=
\lim_{x\downarrow 0}\left(\frac{\overline{\Pi}^{(+)}(x)}
{x^{-1/\kappa} \rho_\kappa(x)+ x^{-2}U_+(x)}\right)
=0.
\end{equation}
It can be proved as in Lemma 4 of \cite{dm2}
that this implies
 \begin{equation}\label{4.5}
\lim_{x\downarrow 0}\left(
\frac{x^{-1/\kappa} \rho_\kappa(x)}
{x^{-2}U_+(x)}
\right)
=0\quad {\rm or}\quad
\liminf_{x\downarrow 0}\left(
\frac{x^{-1/\kappa} \rho_\kappa(x)}
{x^{-2}U_+(x)}
\right)>0.
\end{equation}
Then, since  \eqref{4.3} has been assumed not to hold,
\begin{equation}\label{4.6}
\int_0^{1/2}\frac{x^{\frac{1}{\kappa}}}{\rho_\kappa(x)} \Pi({\rm
d}x)<\infty
\quad {\rm or}\quad
\int_0^{1/2}\frac{x^2}{U_+(x)}\Pi({\rm d}x)<\infty.
 \end{equation}
Noting that, under \eqref{4.1},
 \begin{align*}
\rho_\kappa(x)&=
 \overline {\Pi}^{(+)}(1) -x^{\frac{1}{\kappa}} \overline
{\Pi}^{(+)}(x)
+\int_x^1 y^{\frac{1}{\kappa}}\Pi({\rm d}y)
\nonumber\\
&\le \overline {\Pi}^{(+)}(1) +\int_x^1
y^{\frac{1}{\kappa}}\Pi({\rm d}y)
\sim \int_x^1 y^{\frac{1}{\kappa}}\Pi({\rm d}y), \ {\rm as}\
x\to 0,
\end{align*}
we see that the first relation in \eqref{4.6} is impossible
because it would imply the finiteness of
\[
\int_0^{1/2} x^{\frac{1}{\kappa}}\left(\int_x^1
y^{\frac{1}{\kappa}}\Pi({\rm d}y)\right)^{-1}
\Pi({\rm d}x);
\]
but this is infinite by \eqref{4.1} and the Abel-Dini theorem.
In a similar way, the second relation in \eqref{4.6} can be
shown to be  impossible.
Thus \eqref{4.3} is proved.

Then note that the inverse function $c^\leftarrow$ of $c$ exists
and
satisfies, by \eqref{4.7},
\[
c^\leftarrow(x)
=\frac{\alpha}
{x^{-1/\kappa} \rho_\kappa(x)+ x^{-2}U_+(x)}.
\]
Thus, by \eqref{4.3}, 
\begin{equation}\label{4.7a}
\int_0^1 c^\leftarrow(x) \Pi({\rm d}x)=\infty =
\int_0^{c^\leftarrow(1)}  \overline {\Pi}^{(+)}(c(x)){\rm d}x.
\end{equation}
This proves our claim.\hfill \bbox

\medskip

\begin{prop}\label{pca+} For every $\kappa<1$, \eqref{4.1}
implies
$$\limsup_{t\downarrow 0}\frac{X_t}{t^{\kappa}}=\infty\qquad
\hbox{a.s.}
$$
\end{prop}

\bigskip \noindent{\it Proof of Proposition \ref{pca+}:}\
The  argument relies on the analysis of the completely
asymmetric case when the L\'evy measure $\Pi$ has support in
$[0,1]$ or in $[-1,0]$. Since $\kappa<1$, we can assume
$\gamma=0$ without loss of generality, because of course $\gamma
t=o(t^\kappa)$.  The L\'evy-It\^o decomposition \eqref{eqlid}
then  yields \begin{equation}\label{4.9} X_t=\hat X_t + \tilde
X_t \end{equation} with 
\begin{equation}\label{lid-bis}
\hat X_t=\int_{[0,t]\times[0,1]}xN({\rm d}s,{\rm d}x) \hbox{ and }
\tilde X_t=\int_{[0,t]\times[-1,0]}xN({\rm d}s,{\rm d}x)\,,
\end{equation} where, as usual, the Poissonian stochastic
integrals are taken in the compensated sense.

Choose $\alpha$ so small that
$$(1+\kappa)\alpha<1/2\hbox{ and
}\alpha/(1/2-(1+\kappa)\alpha)^2\le 1/2,$$ and then
$\varepsilon$ so small that $c(\varepsilon)<1$. Observe that for
every $0<t<\varepsilon$  \begin{eqnarray}\label{intens}
t\int_{\{c(t)<x\le 1\}}x\Pi({\rm d}x) &=&tc(t)\overline
{\Pi}^{(+)}(c(t)) +t\lambda(c(t))\nonumber \\ &\le&
\frac{tU_+(c(t))}{c(t)} +t\lambda(c(t))\nonumber \\ & \le&
\alpha(1+\kappa)c(t),  \end{eqnarray} where the last inequality
stems from \eqref{4.7} and \[ \lambda(x):= \int_x^1 \overline
{\Pi}^{(+)}(y){\rm d}y = \int_x^1
y^{1-\frac{1}{\kappa}}y^{\frac{1}{\kappa}-1}
\overline {\Pi}^{(+)}(y){\rm d}y\le \kappa
x^{1-\frac{1}{\kappa}} \rho_\kappa(x) \] (since $\kappa<1$), so
\begin{equation*} \frac{t\lambda(c(t))}{c(t)} \le \frac{\kappa t
\rho_\kappa(c(t))} {c^{\frac{1}{\kappa}}(t)} \le \kappa
\alpha\quad {\rm (by}\ \eqref{4.7}). \end{equation*}

We next deduce
from Lemma \ref{L4.3} that for every $\varepsilon>0$, the
Poisson random measure $N$ has infinitely many atoms in the
domain
$\{(t,x): 0\leq t < \varepsilon\hbox{ and }x>c(t)\}$, a.s.
Introduce $$t_{\varepsilon}:=\sup\{t\leq \varepsilon:
N(\{t\}\times (c(t),1])=1\},$$ the largest instant less than
$\varepsilon$ of such an atom. Our goal is to check that
\begin{equation}\label{eqinq}
P(X_{t_{\varepsilon}-}\geq-c(t_{\varepsilon})/2)\geq 1/33
\end{equation} for every $\varepsilon>0$ sufficiently small, so
that $P(X_{t_{\varepsilon}}>c(t_{\varepsilon})/2)>1/33$.
Since
$t^{\kappa}=o(c(t))$, it follows that
for every  $a>0$
$$P(\exists t\leq \varepsilon: X_t>at^{\kappa})\geq 1/33\,,$$
and hence $\limsup_{t\downarrow 0}  X_t/t^{\kappa}=\infty$ with
probability at least $1/33$. The proof is completed by an appeal
to Blumenthal's $0$-$1$ law.

In order to establish \eqref{eqinq}, we will work henceforth
conditionally on $t_{\varepsilon}$;  recall from the Markov
property of Poisson random measures that the restriction of
$N({\rm d}t,{\rm d}x)$ to $[0,t_{\varepsilon})\times [-1,1]$ is
still a Poisson random measure with intensity ${\rm d}t\Pi({\rm
d}x)$.

Recalling (\ref{intens}) and
discarding  the jumps $\hat \Delta$ of $\hat X$ such that $\hat
\Delta_s>c(t_{\varepsilon})$ for $0\le s< t_{\varepsilon}$ in
the stochastic integral \eqref{lid-bis}, we obtain the
inequality \begin{equation}\label{11z} X_{t_{\varepsilon}-}\geq
\hat
Y_{t_{\varepsilon}-}-\alpha(1+\kappa)c(t_{\varepsilon})+\tilde
X_{t_{\varepsilon}-}
\end{equation}
where $\hat Y_{t_{\varepsilon}-}$ is given by the (compensated)
Poissonian integral $$\hat
Y_{t_{\varepsilon}-}:=\int_{[0,t_{\varepsilon})\times
[0,c(t_{\varepsilon})]}x N({\rm d}s,{\rm
d}x)\,.$$

By a second moment calculation, there is the
inequality
\begin{eqnarray*}
 P\left(
|\hat 
Y_{t_{\varepsilon}-}|>(1/2-\alpha(1+\kappa))c(t_{\varepsilon})\right)&\
\leq
&
\frac{E|\hat
Y_{t_{\varepsilon}-}|^2}{(1/2-\alpha(1+\kappa))^2c^2(t_{\varepsilon})}\\
&\le&  \frac{t_{\varepsilon} \int_{\{0<x\le
c(t_{\varepsilon})\}}x^2\Pi({\rm d}x)}
{(1/2-\alpha(1+\kappa))^2c^2(t_{\varepsilon})}\\
&\le&
\frac{t_{\varepsilon}U_+(c(t_{\varepsilon}))}{(1/2-\alpha(1+\kappa))^2c^2(t_{\varepsilon})}
\\
&\le& \frac{\alpha}{(1/2-\alpha(1+\kappa))^2}\,,\\
\end{eqnarray*}
where the last inequality derives from \eqref{4.7}.
By choice of $\alpha$, the final expression does not exceed
$1/2$. We conclude that \begin{equation}\label{11a}  P\left(
\hat Y_{t_{\varepsilon}-}-\alpha(1+\kappa)c(t_{\varepsilon})\geq
-c(t_{\varepsilon})/2\right)\geq 1/2.
 \end{equation}

We will also use the fact that $\tilde X$
is a mean zero L\'evy process which is
spectrally negative (i.e., with no positive jumps),
so
$$\liminf_{t\downarrow 0}P(\tilde X_t>0) \ge 1/16\,;$$
see \cite{gs}, p. 320. As  furthermore $\tilde X$ is independent
of $\hat Y_{t_{\varepsilon}-}$, we conclude from  \eqref{11z}
and \eqref{11a} that \eqref{eqinq} holds provided that
$\varepsilon$ has been chosen small enough. \hfill\bbox

Now suppose \eqref{4.1} fails.
The remaining results in Theorem \ref{thm3.1}
require the case $\kappa>1/2$ of:

\begin{prop}\label{two}
Assume that $Y$ is a spectrally negative L\'evy process, has
zero mean, and is
not of bounded variation.
Define, for $y>0$, $\lambda >0$,
$$
W_Y(y):=\int_{0}^{y}\int_{x}^{1}z\Pi_Y^{(-)}({\rm d}z){\rm
d}x,$$
and
\begin{equation}\label{JY}
J_Y(\lambda):=\int_0^1\exp \left\{ -\lambda \left(
\frac{y^{{\frac{2\kappa
-1}{\kappa }}}}{W_Y(y)}\right) ^{\frac{\kappa}{1-\kappa }}\right\}
\frac{{\rm d}y}{y},
\end{equation}
where $\Pi_Y^{(-)}$ is the canonical measure of $-Y$,
assumed carried on $(0,1]$, 
and let $\lambda_{Y}^*=\inf\{\lambda>0: J_Y(\lambda)<\infty\}$.
Then with probability one,  for $1/2\le \kappa <1$,
$$\limsup_{t\downarrow 0}\frac{Y_t}{t^\kappa}\ \left\{ \begin{matrix}
=\infty\\
\in(0,\infty)\\
=0\\
\end{matrix}\right. \quad \hbox{according as}\quad 
\lambda_{Y}^*  \left\{ \begin{matrix}
=\infty\\
\in(0,\infty)\\
=0\,.\\
\end{matrix}\right.
 $$
\end{prop}

The proof of Proposition \ref{two} requires several intermediate
steps.
Take  $Y$  as described, then it has characteristic exponent
\[
\Psi_Y(\theta )=\int_{(0,1]}(\e^{-{\rm i}\theta x}-1+ {\rm
i}\theta x)\Pi_Y^{(-)}({\rm d}x).
\]
So we can work with the Laplace exponent
\begin{equation}
\psi_Y(\theta )=\Psi_Y(-{\rm i}\theta)=\int_{(0,1]}(\e^{-\theta
x}-1+\theta x)\Pi_Y^{(-)}({\rm d}x),  \label{x1}
\end{equation}
such that $E\e^{\theta Y_t}=\e^{t\psi_Y(\theta)}$, $t\ge 0$,
$\theta\ge 0$.

Let $T=(T_t,t\geq 0)$ denote the first
passage
process of $Y;$ this is a subordinator whose Laplace exponent
$\Phi $ is
the inverse function to $\psi_Y$
(\cite{b}, p.189), and since $Y(T_{t})\equiv t$ we see that
the alternatives in Proposition \ref{two}
can be deduced immediately from
$$\limsup_{t\downarrow 0}\frac{Y_t}{t^\kappa}\ \left\{ \begin{matrix}
=\infty\\
\in(0,\infty)\\
=0\\
\end{matrix}\right. \ \Longleftrightarrow\ 
\liminf_{t\downarrow 0}\frac{T_{t}}{t^{1/\kappa }}\left\{
\begin{matrix}
=0\\
\in(0,\infty)\\
=\infty.\\
\end{matrix}\right. $$

The subordinator
$T$ must have zero drift since if
$\lim_{t\downarrow 0}T_t/t:=c>0$ a.s. then
$\sup_{0<s\le T_t}Y_s =t$ (see \cite{b}, p.191) would give
$\limsup_{t\downarrow 0}Y_t/t\le 1/c<\infty$ a.s., thus $Y\in bv$,
which is
not the case. We can assume $T$ has no
 jumps bigger  than 1,
and further exclude the trivial case when
$T$ is compound Poisson. So the main part of
the proof of Proposition \ref{two} is
the following, which is a kind of analogue of Theorem 1 of 
Zhang \cite{zhang}.

\begin{lem}\label{ronlem1}
Let $T$ be any subordinator with zero drift whose L\'{e}vy
measure $\Pi_T$ is
carried by $(0,1]$ and has $\overline\Pi_T(0+)=\infty$,
where $\overline\Pi_T(x)=\Pi_T\{(x,\infty)\}$ for $x>0$.
Put $m_T(x)=\int_{0}^{x}\overline\Pi_T(y){\rm d}y$ 
and for $d>0$ let
\begin{equation}\label{KT}
K_T(d):=\int_{0}^{1}\frac{{\rm d}y}{y}\exp \left\{
-d\frac{\left(m_T(y)\right)^{\frac{\gamma}{\gamma-1}}}{y}\right\},\
{\rm where}\ \gamma>1.
\end{equation}
Let
$d_K^*:= \inf\{d>0: K_T(d)<\infty\}\in[0,\infty]$.
Then, with probability one,
\newline (i)\ $d_K^*=0$ iff 
$\displaystyle{\lim_{t\downarrow 0}
\frac{T_t}{t^\gamma}=\infty}$;
\newline (ii)\ $d_K^*=\infty$ iff 
$\displaystyle{\liminf_{t\downarrow 0} \frac{T_t}{t^\gamma}=0}$;
\newline (iii)\ $d_K^*\in (0,\infty)$ iff 
$\displaystyle{\liminf_{t\downarrow 0} \frac{T_t}{t^\gamma}=c}$,
for some $c\in(0,\infty)$.
\end{lem}

Before beginning the proof of Lemma \ref{ronlem1}, we need
some preliminary results.
To start with, we need the following lemma.

\begin{lem}\label{ronlem2}
Let $S=(S_t, t\geq 0)$ be a subordinator,
and $a$ and $\gamma$ positive constants.
Then
\begin{equation}\label{x22}
\liminf_{t\downarrow 0}\frac{S_t}{t^\gamma}\le a\ {\rm a.s.}
\end{equation}
if and only if for every $r\in(0,1)$ and $\eta>0$
\[
\sum_{n\ge 1}
P(S_{r^n}\leq (a+\eta)r^{n\gamma }\text{ })=\infty.\]
\end{lem}

\bigskip \noindent{\it Proof of Lemma \ref{ronlem2}:}\
One way is obvious, so suppose 
\begin{equation}\label{rev0}
\sum P\left(S_{r^n}\leq ar^{n\gamma}/(1-r)\right)=\infty.
\end{equation}
For a given $\veps>0$, define events
\[
A_n=\{S_{r^n}-S_{r^{n+1}}\le ar^{n\gamma}/(1-r)\},\
B_n=\{S_{r^{n+1}}\le \veps r^{(n+1)\gamma}\},\
n\ge 0.
\]
Then the $\{A_n\}_{n\ge 0}$ are independent, and each $B_n$ is
independent of
the collection 
$\{A_n, A_{n-1}, \cdots, A_0\}$.
Further, $\sum_{n\ge 0}P(A_n)=\infty$ by \eqref{rev0} 
(recall $S$ is a subordinator), so
$P(A_n\ {\rm i.o.})=1$.
Then, by the Feller--Chung lemma (\cite{CT}, p. 69)
we can deduce that $P(A_n\cap B_n\ {\rm i.o.})=1$, provided
$P(B_n)$ 
is bounded away from 0: $P(B_n)\ge 1/2$, say, for $n$ large
enough.
To see that this is the case here, take $b>0$ and
$\veps\in(0,1)$
and truncate the jumps of $S$ (which is of bounded variation) at
$b\veps>0$,
where $b$ will be specified more precisely shortly.
Thus, let
$S_t^\veps
=\sum_{0<s\le t} \Delta S_s \boldsymbol{1}_{\{\Delta S_s \le
b\veps)\}}$.
Now $S_t-S_t^\veps$ is nonzero only if there is at least one
jump in  $S$ of magnitude greater than $b\veps$ 
up till time $t$, and this has probability
bounded above by $t\overline \Pi_S(b\veps)$, 
where $\Pi_S$ is the L\'evy measure of
$S$.
Then by a standard truncation argument, and using a first-moment
Markov
inequality,
 \begin{eqnarray*}
P(S_t>\veps t^\gamma)
&\le& 
\frac{t\int_{(0,b\veps]}x\Pi_S({\rm d}x)}{\veps t^\gamma}
+t\overline
\Pi_S(b\veps)
\nonumber \\
&\le&
\frac{t\int_{(0,b\veps]}x\Pi_S({\rm d}x)}{\veps t^\gamma} +
\frac{t\overline \Pi_S(b\veps)}{\veps t^\gamma}
\quad {\rm (once}\ \veps t^\gamma\le 1)
\nonumber \\
&=&
\frac{tm_T(b\veps)}{\veps t^\gamma}
\le \frac{tm_T(b)}{\veps t^\gamma}.
\end{eqnarray*}
Now choose $b=h(\veps t^{\gamma-1}/2)$, where $h(\cdot)$ is the
inverse function
to $m_T(\cdot)$.
Then the last ratio is smaller than $1/2$.
Replacing $t$ by $r^{n+1}$ in this gives
$P(B_n)\ge 1/2$ for $n$ large enough.
Finally $P(A_n\cap B_n\ {\rm i.o.})=1$ implies
$P(S_{r^n}\le r^{n\gamma}(a/(1-r)+\veps r^\gamma)\ {\rm i.o.})=1$,
thus
\[
\liminf_{t\downarrow 0}
\frac{S_t}{t^\gamma} \le \frac{a}{1-r}+\veps r^\gamma\ {\rm a.s.},
\]
in which we can let $\veps\downarrow 0$ and $r\downarrow 0$ to
get
\eqref{x22}.
\hfill\bbox

Applying Lemma \ref{ronlem2} to $T_t$, we see that the
alternatives in (i)--(iii) of Lemma \ref{ronlem1}
hold iff for some $r<1$,
for all, none, or some but not all, $a>0$,
\begin{equation}
\sum_{n\ge 1} P(T_{r^{n}}\leq ar^{n \gamma }\text{ })<\infty.
\label{x3}
\end{equation}%

The next step is to get bounds for the probability in
\eqref{x3}. One
way is easy. Since $\overline\Pi_T(0+)=\infty$,
$\overline\Pi_T(x)$ 
is strictly positive, and thus $m_T(x)$ is strictly
increasing, on a neighbourhood of 0.
Recall that we write $h(\cdot)$ for the inverse function to
$m_T(\cdot)$.

\begin{lem}\label{ronlem3}
Let $T$ be a subordinator with canonical measure $\Pi_T$
satisfying $\overline\Pi_T(0+)=\infty$.
Then there is an absolute constant $K$ such that,
for any $ c>0$ and $\gamma>0$,
\begin{equation}
P(T_{t}\leq ct^{\gamma })\leq 
\exp \left\{-\frac{ct^{\gamma }}{h(2ct^{\gamma -1}/K)}\right\},
\text{ }t>0.  \label{x4}
\end{equation}
\end{lem}

\bigskip \noindent{\it Proof of Lemma \ref{ronlem3}:}\
We can write
\[
\Phi(\lambda)= -\frac{1}{t}\log E\e^{-\lambda T_t}
=\int_{(0,1]} \left(1-\e^{-\lambda x}\right)\Pi_T({\rm d}x),\
\lambda>0.
\]
Markov's inequality gives, for any $\lambda >0,$
$c>0$,
\begin{eqnarray*}
P(T_{t} \leq ct^{\gamma })&\leq& \e^{\lambda ct^{\gamma
}}E(\e^{-\lambda T_{t}})\\
&\leq&
\exp \{-\lambda t(\lambda ^{-1}\Phi(\lambda )-ct^{\gamma-1})\} 
\\
&\leq &
\exp \{-\lambda t(Km_T(1/\lambda )-ct^{\gamma-1})\},
\ {\rm for\ some}\ K>0,
\end{eqnarray*}
where we have used \cite{b}, Prop. 1, p. 74. 
Now choose $\lambda =1/h(2ct^{\gamma -1}/K)$ and we have
(\ref{x4}).
\hfill\bbox

The corresponding lower bound is trickier:

\begin{lem}\label{ronlem4}
Suppose that $T$ is as in Lemma \ref{ronlem1}, and additionally
satisfies $%
\lim_{t\downarrow 0}P(T_{t}\leq dt^{\gamma })=0$ for some $d>0$
and $\gamma>1$.
Then for
any $c>0$%
\begin{equation}
P(T_{t}\leq ct^{\gamma })\geq \frac{1}{4}
\exp\left\{-\frac{ct^{\gamma }}{h(ct^{\gamma -1}/4)}\right\}
\text{ for all small enough }t>0.
\label{x5}
\end{equation}
\end{lem}

\bigskip \noindent{\it Proof of Lemma \ref{ronlem4}:}\
Take $\gamma>1$ and assume $\lim_{t\downarrow 0}P(T_{t}\leq
dt^{\gamma })=0$, where $d>0$.
First we show that 
\begin{equation}
\frac{t^{\gamma }}{h(t^{\gamma -1})}\rightarrow \infty \text{ as
}%
t\downarrow 0.  \label{x6}
\end{equation}%
To do this we write, for each fixed $t>0$,
$T_{t}=T_{t}^{(1)}+T_{t}^{(2)}$, where the distributions of the
independent random variables $T_{t}^{(1)}$ and $T_{t}^{(2)}$ are
specified by
\[
\log E\e^{-\lambda T_{t}^{(1)}} =-t\int_{(0, h(\varepsilon
t^{\gamma-1})]}(1-\e^{-\lambda x})\Pi_T({\rm d}x)
\]
and
\[
\log E\e^{-\lambda T_{t}^{(2)}} =-t\int_{(h(\varepsilon
t^{\gamma -1}),1]}(1-\e^{-\lambda x})\Pi_T({\rm d}x),
\]
for a given $\varepsilon >0$.
Observe that 
$$ET_{t}^{(1)}=t\int_{(0, h(\varepsilon
t^{\gamma
-1})]}x\Pi_T
({\rm d}x)\leq tm_T(h(\varepsilon t^{\gamma -1}))=\varepsilon
t^{\gamma
},$$
so that 
\begin{eqnarray*}
P(T_{t} >dt^{\gamma })&\leq& P(T_{t}^{(1)}>dt^{\gamma
})+P(T_{t}^{(2)}\neq 0)
\\
&\leq&\frac{ET_{t}^{(1)}}{dt^{\gamma }}+1-P(T_{t}^{(2)}=0)\\
&\leq&
\varepsilon
/d+1-P(T_{t}^{(2)}=0).
\end{eqnarray*}%
Thus for all sufficiently small t,
\[
P(T_{t}^{(2)}=0)\leq \varepsilon /d+P(T_{t}\leq dt^{\gamma
})\leq
2\varepsilon /d,\]
because of our assumption that
$\lim_{t\downarrow 0}P(T_{t}\leq dt^{\gamma })=0$.
Now 
$$P(T_{t}^{(2)}=0)=\exp\left(-t\overline\Pi_T(h(\varepsilon
t^{\gamma -1}))\right)$$
and 
\[
t\overline\Pi_T(h(\varepsilon t^{\gamma -1}))\leq
\frac{t}{h(\varepsilon
t^{\gamma -1})}m_T(h(\varepsilon t^{\gamma
-1}))=\frac{\varepsilon
t^{\gamma }%
}{h(\varepsilon t^{\gamma -1})},
\]%
so we see, taking say $\varepsilon =d/4,$ that $h(t^{\gamma
-1})\leq at^{\gamma}$ for a constant $a>0$,
or, equivalently, $h(t)\leq at^{1+\beta}$, where
$\beta=1/(\gamma -1)>0$, for all sufficiently small $t$.
However, $m_T(\cdot)$ is
concave, so its inverse function $h$ is convex, so $h(t/2)\leq
h(t)/2\leq
a(1/2)^{\beta }(t/2)^{1+\beta }$,
or
$h(t)\le a(1/2)^{\beta }t^{1+\beta }$,
for small $t.$ Iterating this argument gives (\ref{x6}).

Now write $\eta =\eta (t)=h(ct^{\gamma -1}/4)$ and
define processes $(Y_t^{(i)})_{t\ge 0}$, $i=1,2,3$, such
that $(T_t)_{t\ge 0}$ and $(Y_t^{(1)})_{t\ge 0}$ are independent,
and
$(Y_t^{(2)})_{t\ge 0}$ and $(Y_t^{(3)})_{t\ge 0}$ are independent,
and are such that
$\log E(\e^{-\lambda Y_{t}^{(i)}})=-t\int_{(0,1]}
(1-\e^{-\lambda x})\Pi_T
^{(i)}({\rm d}x),$ $%
i=1,2,3,$ where%
\begin{eqnarray*}
\Pi_T^{(1)}({\rm d}x) &=&\overline\Pi_T(\eta )\delta_{\eta
}({\rm d}x), \\
\Pi_T^{(2)}({\rm d}x) &=&\overline\Pi_T(\eta )\delta_{\eta
}({\rm d}x)+\boldsymbol{1}%
_{(0,\eta ]}\Pi_T({\rm d}x), \\
\Pi_T^{(3)}({\rm d}x) &=&\boldsymbol{1}_{(\eta, 1]}\Pi_T({\rm
d}x),
\end{eqnarray*}%
and $\delta_{\eta }({\rm d}x)$ is the point mass at $\eta$.
Then we have
$T_{t}+Y_{t}^{(1)}\overset{d}{=}Y_{t}^{(2)}+Y_{t}^{(3)}$, and 
 \begin{eqnarray*}
P(T_{t} \leq ct^{\gamma })
&\geq&
P(T_{t}+Y_{t}^{(1)} \leq ct^{\gamma })
\\
&\geq&
P(Y_{t}^{(3)}=0)P(Y_{t}^{(2)}\leq ct^{\gamma })\\
&=&\e^{-t\overline\Pi_T(\eta )}P(Y_{t}^{(2)}\leq ct^{\gamma }).
 \end{eqnarray*}
Since $t\overline\Pi_T(\eta )=\eta ^{-1}t\eta
\overline\Pi_T(\eta )\leq
\eta ^{-1}tm_T(\eta )\le ct^{\gamma }/h(ct^{\gamma-1}/4),$
(\ref{x5})
will follow when we show that $\liminf_{t\downarrow 0}P(Y_{t}^{(2)}\leq
ct^{\gamma })\geq
1/4.$ By construction we have 
$$
EY_{t}^{(2)}=t\left(\int_{0}^{\eta }x\Pi_T({\rm d}x)+\eta
\overline\Pi_T(\eta )\right) 
=tm_T(\eta )=\frac{ct^{\gamma }}{4},
$$
so if we put $Z_{t}=Y_{t}^{(2)}-ct^{\gamma }/4$ and
write $t\sigma_t^2=EZ_t^2$ we can apply Chebychev to get 
\[
P\left(Y_{t}^{(2)}\leq ct^{\gamma }\right)\geq 
P\left(Z_t\le \frac{ct^{\gamma }}{2}\right)
\geq \frac{5}{9}
\]
when $t\sigma _{t}^{2}\leq c^{2}t^{2\gamma }/9.$ 
To deal with the opposite case,  $t\sigma _{t}^{2}>
c^{2}t^{2\gamma }/9$,
 we use  the Normal approximation in Lemma \ref{ron2lem1}. 
In the notation of that lemma, 
$m_3=\int_{|x|\le \eta}|x|^3 \Pi_T^{(2)}({\rm d}x)$
and
$m_2=\sigma^2_t$, in the present situation.
Since, then, $m_3\le \eta \sigma^2_t$, and
we have 
$\eta=h(ct^{\gamma-1}/4)=o(t^\gamma)$, as $t\downarrow 0$,
by \eqref{x6},
we get 
\[
\sup_{x\in \R}
\left|P(Z_{t}\leq x\sqrt{t}\sigma _t)-1+\overline{F}(x)\right|
\leq 
\frac{A\eta \sigma^2_t}{\sqrt{t} \sigma^3_t} =o(1),\
{\rm as}\ t\downarrow 0.
\]
Choosing $x=ct^{\gamma }/(2\sqrt{t}\sigma _{t})$ gives
$P(Z_{t}\leq
ct^{\gamma }/2)\geq 1/4$, hence \eqref{x5}.  \hfill\bbox

\medskip

We are now able to establish Lemma \ref{ronlem1}.
Again, recall, $h(\cdot)$ is inverse to $m_T(\cdot)$.

\bigskip \noindent{\it Proof of Lemma \ref{ronlem1}:}\
(i)\
Suppose that $K_T(d)<\infty $ for some $d>0$
and write $x_{n}=h(cr^{n(\gamma -1)})$, where $\gamma>1$ and
$0<r<1$.
Note that since $h(x)/x$ is increasing we have $%
x_{n+1}\leq r^{\gamma -1}x_{n}.$ Also we have
$m_T(x_{n})=Rm_T(x_{n+1})$ where $%
R=r^{1-\gamma }>1.$ So for $y\in \lbrack x_{n+1},x_{n}),$
\begin{eqnarray*}
\frac{m_T(y)^{\frac{\gamma }{\gamma -1}}}{y} &\leq
&\frac{m_T(x_{n})^{^{\frac{%
\gamma }{\gamma -1}}}}{x_{n+1}}\\
&=&\frac{R^{\frac{\gamma }{\gamma -1}%
}m_T(x_{n+1})^{^{\frac{\gamma }{\gamma -1}}}}{x_{n+1}} \\
&=&\frac{(Rc)^{^{^{\frac{\gamma }{\gamma -1}}}}r^{(n+1)\gamma
}}{%
h(cr^{(n+1)(\gamma -1)})}.
\end{eqnarray*}%
Thus
\begin{eqnarray*}
\int_{x_{n+1}}^{x_{n}}\frac{{\rm d}y}{y}\exp \left\{
-d\frac{m_T(y)^{\frac{\gamma}{\gamma-1}}   }{y}%
\right\} 
&\geq &
\exp \left\{ -\frac{d(Rc)^{^{^{\frac{\gamma }{\gamma -1}%
}}}r^{(n+1)\gamma }}{h(c r^{(n+1)(\gamma -1)})}\right\} \log
\frac{x_{n}}{%
x_{n+1}} \\
&\ge&(\log R)\exp \left\{ -\frac{c^{\prime }r^{(n+1)\gamma
}}{h(2c^{\prime
}r^{(n+1)(\gamma -1)}/K)}\right\},
\end{eqnarray*}%
where $K$ is the constant in Lemma \ref{ronlem3}
and we have chosen%
\begin{equation*}
c=\left( \frac{K}{2d}\right) ^{\gamma -1}R^{-\gamma }\text{ and
}c^{\prime
}=Kc/2.
\end{equation*}
Then $K_T(d)<\infty$ gives
$\sum_{1}^{\infty }P(T_{r^{n}}\leq c^{\prime }r^{n\gamma
})<\infty,$
and so 
$\liminf_{t\downarrow 0}T_{t}/t^\gamma \geq c^{\prime }>0$ a.s.,
by Lemma \ref{ronlem2}.
Thus we see that $\liminf_{t\downarrow 0}T_{t}/t^\gamma =0$ a.s.
implies that $K_T(d)=\infty $ for every $d>0$, hence
$d_K^*=\infty$.

Conversely, assume that 
$\liminf_{t\downarrow 0}T_{t}/t^\gamma >0$ a.s.
Then by Lemma \ref{ronlem2},
$\sum_{1}^{\infty }P(T_{r^{n}}\leq
cr^{n\gamma })<\infty $ for some $c>0$ and $0<r<1$.
Then $P(T_{r^{n}}\leq
cr^{n\gamma })\rightarrow 0,$ Lemma \ref{ronlem4} applies,
 and we have
\[
\sum_{1}^{\infty }\exp\left\{-\frac{cr^{n\gamma
}}{h(cr^{n(\gamma
-1)}/4)}\right\}%
<\infty.
\]%
Putting $x_{n}=h(cr^{n(\gamma -1)}/4)$
(similar to but not the same $x_n$ as in the previous paragraph),
and $c^\prime=4^{\frac{\gamma }{\gamma -1}}/c^{\frac{1}{\gamma -1}}$,
we see that
\begin{equation}
\sum_{1}^{\infty }\exp \left\{
-\frac{c^{\prime }m_T(x_{n})^{\frac{\gamma
}{\gamma -1}}}{x_{n}}\right\}<\infty.  \label{x7}
\end{equation}%
We have $m_T(x_{n-1})=Rm_T(x_{n})$
where $R=r^{1-\gamma }>1.$ Take $L>R$ and let $k_{n}=\min
(k\ge 1:x_{n-1}L^{-k}\leq
x_{n}),$ so that $x_{n-1}L^{-k_{n}}\le x_{n}$. 
Then for any $d>0$
\begin{eqnarray*}
&&\int_{x_{n}}^{x_{n-1}}\exp\left\{
-d\frac{m_T(y)^{\frac{\gamma }{\gamma -1}}}{y}\right\}
y^{-1}{\rm d}y
\\
&\leq&\sum_{i=1}^{k_{n}}\int_{x_{n-1}L^{-i}}^{x_{n-1}L^{1-i}}
\exp\left\{
-dm_T(y)^{\frac{1}{\gamma-1}}\frac{m_T(y)}{y}\right\}y^{-1}{\rm
d}y
\\
&\leq&
\sum_{i=1}^{k_{n}}\int_{x_{n-1}L^{-i}}^{x_{n-1}L^{1-i}}\exp
\left\{-dm_T(x_{n})^{\frac{1}{\gamma-1}}
\frac{m_T(x_{n-1}L^{-i})}{x_{n-1}L^{-i}} \right\}y^{-1}{\rm d}y
\\
&\le&
\log L\sum_{i=1}^{\infty }\exp\left\{
-dm_T(x_{n})^{\frac{1}{\gamma
-1}}\frac{L^im_T(x_{n}L^{-1})}{x_{n-1}} 
\right\}
\\
&\leq&
\log L\sum_{i=1}^{\infty }\exp \left\{
-dm_T(x_{n})^{\frac{1}{\gamma-1}}\frac{L^{i-1}m_T(x_{n})}{x_{n-1}} 
\right\}
\\
&=&
\log L\sum_{i=1}^{\infty }\exp \left\{ -dL^{i-1}R^{-\frac{\gamma
}{\gamma-1}}
\frac{m_T(x_{n-1})^{\frac{\gamma }{\gamma -1}}}{x_{n-1}}\right\}.
\end{eqnarray*}
Approximate this last sum with an integral of the form
$\int_0^\infty a_n^{L^x}{\rm d}x$, where 
$a_n=\exp(-c^\prime m_T(x_{n-1})^{\frac{\gamma}{\gamma -1}}/x_{n-1})$,
with $d=c^\prime LR^{\frac{\gamma }{\gamma-1}}
=4^{\frac{\gamma }{\gamma-1}}Lr^{-\gamma}/c^{\frac{1}{\gamma-1}}$, 
to see that it is bounded above by a constant
multiple of $a_n$. It follows from (\ref{x7}) that
$\sum a_n<\infty$, hence we get
$K_T(d)<\infty $, and Part (i) follows.

(ii)\ If $K_T(d)<\infty$ for all $d>0$ then, because
$c^{\prime }\rightarrow 0$ as $d\rightarrow \infty$
at the end of the proof of the forward part of Part (i), we have
$\liminf_{t\downarrow 0}T_t/t^\gamma=\infty $ a.s., i.e.,
$\lim_{t\downarrow 0}T_t/t^\gamma=\infty $ a.s.
Conversely, if this holds, then because
$d\rightarrow \infty$ as $c\rightarrow 0$
at the end of the proof of the converse part of
Part (i), we have $K_T(d)<\infty$ for all $d>0$.
This completes the proof of Lemma \ref{ronlem1}.
\hfill\bbox

\bigskip \noindent{\it Proof of Proposition \ref{two}:}\
To finish the proof of Proposition \ref{two}, we need only show
that
$J_Y(\lambda)=\infty $ for all $\lambda>0$ is equivalent to 
$K_T(d)=\infty $ for all $d>0$,
where $K_T(d)$ is evaluated for the first-passage process $T$
of $Y,$ and $ \gamma =1/\kappa$. 
We have  from \eqref{x1}, after integrating by
parts,
\begin{eqnarray*}
\psi_Y(\theta ) &=&\int_{0}^{1}(\e^{-\theta x}-1+\theta
x)\Pi_Y^{(-)}({\rm d}x) 
\nonumber \\
&=&\theta
\int_{0}^{1}(1-\e^{-\theta x})\overline\Pi_Y^{(-)}(x){\rm d}x,
\end{eqnarray*}
and differentiating \eqref{x1} gives
\[
\psi_Y'(\theta) 
= \int_0^1x(1- \e^{-\theta x})\Pi_Y^{(-)}({\rm d}x).
 \]
So we see that $\theta ^{-1}\psi _{Y}(\theta )$ and $\psi
_{Y}^{\prime
}(\theta )$ are Laplace exponents of driftless subordinators,
and using the
estimate in \cite{b}, p.74, twice, we get 
\begin{equation*}
\psi _{Y}(\theta )\asymp \theta ^{2}\widetilde{W}_{Y}(1/\theta
)\text{ and }
\psi _{Y}^{\prime }(\theta )\asymp \theta W_{Y}(1/\theta ),
\end{equation*}
where  $\widetilde{W_{Y}}(x)=\int_{0}^{x}A_{Y}(y)\mathrm{d}y$ and
$A_{Y}(x):=\int_{x}^{1}\overline{\Pi }_{Y}^{(-)}(y)\mathrm{d}y$,
for $x>0$,
we recall the definition of $W_{Y}$ just prior to \eqref{JY},
and $``\asymp "$ means that the ratio of the quantities on each side 
of the symbol is bounded above and
below by finite positive constants for all values of the
argument. However,
putting $U_{Y}(x)=\int_{0}^{x}2z\overline{\Pi
}_{Y}^{(-)}(z)\mathrm{d}z$, for $x>0$, 
we see that 
\[
W_Y(x)=\int_0^x \int_y^1z \Pi_Y^{(-)}({\rm d}z){\rm d}y
 =\frac{1}{2}U_Y(x)+\tilde W_Y(x),
\]
and 
\[
\tilde W_Y(x)=\frac{1}{2} U_Y(x)+xA_Y(x);
\]
thus 
\[
\tilde W_Y(x)\le
W_Y(x) =U_Y(x)+xA_Y(x)\le 2\tilde W_Y(x).
\]
Hence we have
\begin{equation}
\theta ^{2}W_{Y}(1/\theta )\asymp \psi _{Y}(\theta )\asymp
\theta \psi_{Y}^{\prime }(\theta ).  \label{inW}
\end{equation}
We deduce that $J_Y(\lambda)=\infty $ for all $\lambda >0$ is
equivalent to $ \tilde J_Y(\lambda)=\infty $ for all
$\lambda>0$,
where
\begin{eqnarray*}
\tilde J_Y(\lambda) &=&\int_{0}^{1}\exp \left\{ -\lambda
y^{\frac{-1}{1-\kappa 
}}\psi_Y (1/y)^{\frac{-\kappa }{1-\kappa }}\right\} \frac{{\rm
d}y}{y}
\\
&=&\int_{1}^{\infty }\exp \left\{ -\lambda y^{\frac{1}{1-\kappa
}}\psi (y)^{%
\frac{-\kappa }{1-\kappa }}\right\} \frac{{\rm d}y}{y}.
\end{eqnarray*}%
But we know that $\Phi$, the exponent of the first-passage
process $T$, is
the inverse of $\psi_Y$, so making the obvious change of
variable gives%
\[
\tilde J_Y(\lambda)=\int_{\psi_Y (1)}^{\infty }\exp \left\{
-\lambda
\Phi (z)^{%
\frac{1}{1-\kappa }}z^{\frac{-\kappa }{1-\kappa }}\right\}
\frac{\Phi ^{\prime }(z){\rm d}z}{\Phi (z)}.
\]
{}From \eqref{inW} we deduce that $z\Phi ^{\prime}(z)/\Phi(z)\asymp
1$ for all $z>0$, so $J_Y(\lambda)=\infty $ for all $\lambda>0$ is
equivalent to $\hat J_Y(\lambda)=\infty $ for all 
$\lambda>0,$ where
\begin{eqnarray*}
\hat J_Y(\lambda) &=&\int_{1}^{\infty }\exp \left\{ -\lambda
\Phi (z)^{\frac{%
1}{1-\kappa }}z^{\frac{-\kappa }{1-\kappa }}\right\}
\frac{{\rm d}z}{z}, \\
&=&\int_{0}^{1}\exp \left\{ -\lambda \Phi
(z^{-1})^{\frac{1}{1-\kappa }}z^{%
\frac{\kappa }{1-\kappa }}\right\} \frac{{\rm d}z}{z}.
\end{eqnarray*}%
Since $\Phi (z^{-1})$ is bounded above and below by multiples of
$z^{-1}m_T(z)$ (\cite{b}, p.74), our claim is established.
\hfill\bbox

We are now able to complete the proof of Theorem \ref{thm3.1}:

\bigskip \noindent{\it Proof of Theorem \ref{thm3.1}:}\
The implications (i)$\Rightarrow$\eqref{1s} and
(iii)$\Rightarrow$\eqref{ze}
 stem from
Proposition \ref{pca+} and Theorem \ref{th1}, respectively.
So we can focus on the situation when
 $$\int_0^1 \pibar^{(+)}(x^\kappa){\rm d}x<\infty=\int_0^1
\pibar^{(-)}(x^\kappa){\rm d}x\,.$$
Recall the decomposition \eqref{4.9}
where $\hat X_t$ has canonical measure $\Pi^{(+)}({\rm d}x)$.
Thus from Theorem \ref{th1}, $\hat X_t$ is $o(t^\kappa)$ a.s. 
as $t\downarrow 0$.
Further, $\tilde X$ is spectrally negative with mean zero.
When (ii), (iv) or (v) holds, $X\notin bv$ (see Remark
\ref{r40} (ii)), so $X_t/t$ takes arbitrarily large
positive and negative values a.s. as $t\downarrow 0$, and
$\liminf_{t\downarrow 0}X_t/t^\kappa \le 0 \le 
\limsup_{t\downarrow 0}X_t/t^\kappa$ a.s.
The implications (ii)$\Rightarrow$\eqref{1s},
(iv)$\Rightarrow$\eqref{ze}
and (v)$\Rightarrow$\eqref{ze2}
now follow from Proposition \ref{two}.
\hfill\bbox

\begin{remark}\label{wh?}\ 
Concerning Remark \ref{r40} (iii):
perusal of the proof of Theorem \ref{th2} shows that we can add
to $X$ a compound Poisson process with masses $f_\pm(t)$, say,
at $\pm\sqrt{t}$, provided 
$\sum_{n\ge 1} \sqrt{t_n}f_\pm(t_n)$ converges, 
and the proof remains valid.
The effect of this is essentially only to change the kind of
truncation that is being applied, without changing the value of
the limsup,
and in the final result this
shows up only in an alteration to $V(x)$.
Choosing $f(t)$ appropriately, the new $V(\cdot)$ becomes
$U(\cdot)$ or $W(\cdot)$, which are thus equivalent in the
context of Theorem \ref{th2}.
Note that we allow $\kappa=1/2$ in Proposition \ref{two}.
We will omit further details,
but the above shows there is no contradiction with Theorem \ref{th2}.

\end{remark}

\subsection{Proof of Theorem \ref{thm3.2}}\
Theorem \ref{thm3.2} follows by taking $a(x)=x^\kappa$,
$\kappa>1$,  in the
Propositions \ref{prop3} and \ref{prop3-end} below, which are a
kind of generalisation of
Theorem 9 in Ch. III of \cite{b}. Recall the definition of
$A_-(\cdot)$ in \eqref{A-def}.\hfill\bbox

\begin{prop}\label{prop3}
Assume $X\in bv$ and $\delta=0$.
Suppose $a(x)$ is a positive
deterministic measurable function on $[0,\infty)$ with
$a(x)/x$ nondecreasing and $a(0)=0$.
Let $a^\leftarrow(x)$ be its inverse function.
Suppose 
\begin{equation}\label{a2}
\int_0^1 \frac{\Pi({\rm d}x)}{ 1/a^\leftarrow(x)+A_-(x)/x}
=\infty.
\end{equation}
Then we have
\begin{equation}\label{a1}
\limsup_{t\downarrow 0}
\frac{X_t}{a(t)}=\ \infty\ {\rm a.s.}
\end{equation}
\end{prop}

\bigskip \noindent{\it Proof of Proposition \ref{prop3}:}\quad
Assume $X$ and $a$ as specified.
Then the function  $a(x)$ is strictly increasing, so $a^\leftarrow(x)$
is well
defined, positive, continuous, and nondecreasing on
$[0,\infty)$, with $a^\leftarrow(0)=0$ and 
$a^\leftarrow(\infty)=\infty$.
Note that the function
\[
\frac{1}{a^\leftarrow(x)} +\frac{A_-(x)}{x}
= \frac{1}{a^\leftarrow(x)}+\int_0^1\overline{\Pi}^{(-)}(xy){\rm
d}y
\]
is continuous and nonincreasing, tends to $\infty$ as $x\to
0$, and to 0 as $x\to \infty$.
Choose  $\alpha\in(0,1/2)$ arbitrarily small so that
$$2 \alpha\left(2(1/2-\alpha)^{-2}+1\right) \le 1,$$ 
and define, for $t>0$,
\[
b(t)=\inf\left\{x>0: 
\frac{1}{a^\leftarrow(x)} +\frac{A_-(x)}{x}
\le \frac{\alpha}{t}\right\}.
\]
Then $0<b(t)<\infty$ for $t>0$, $b(t)$ is strictly increasing,
$\lim_{t\downarrow 0}b(t)=0$,  and
\begin{equation}\label{5.11}
\frac{t}{a^\leftarrow(b(t))}  +\frac{tA_-(b(t))}{b(t)}= \alpha.
\end{equation}
Also
 $b(t)\ge a(t/\alpha)$, and 
the inverse function $b^\leftarrow(x)$ exists and
satisfies
 \begin{equation}\label{bdef}
b^\leftarrow(x)
=\frac{\alpha}
{1/a^\leftarrow(x)+ A_-(x)/x}.
 \end{equation}
Thus, by \eqref{a2},
\begin{equation}\label{5.2}
\int_0^1 b^\leftarrow(x) \Pi({\rm d}x)=\infty =
\int_0^1 \overline {\Pi}^{(+)}(b(x)){\rm d}x.
\end{equation}

Set
$$U_-(x):=2\int_0^x y\overline{\Pi}^{(-)}(y) {\rm d}y.$$
Then we have the upper-bounds
\[
t\overline{\Pi}^{(-)}(b(t))
\le \frac{tA_-(b(t))}{b(t)} \le \alpha,
\]
and
\[
\frac{tU_-(b(t))}{b^2(t)}
\le \frac{2tA_-(b(t))}{b(t)}
\le 2\alpha.
\]

Since $X\in bv$ and $\delta=0$  we can express
$X$ in terms of its positive and negative jumps,
$\Delta_s^{(+)}=\max(0,\Delta_s)$ and
$\Delta_s^{(-)}=\Delta_s^{(+)}-\Delta_s$:
\begin{equation}\label{5.3}
X_t= \sum_{0<s\le t} \Delta_s^{(+)}
-\sum_{0<s\le t} \Delta_s^{(-)} 
= X_t^{(+)} - X_t^{(-)}, \ {\rm say}.
\end{equation}
Recall that $\Delta_s^{(\pm)} \le 1$ a.s.
We then have
\begin{eqnarray*}
& &P(X_t^{(-)}  > b(t)/2) \\
&\le & P\left(\sum_{0<s\le t} (\Delta_s^{(-)}\wedge b(t))
>b(t)/2\right)
+ P\left(\Delta_s^{(-)}>b(t)\ {\rm for\ some}\ s\le t \right)
 \\
&\le& P\left(\sum_{0<s\le t}(\Delta_s^{(-)}\wedge
b(t))-tA_-(b(t))>(1/2-\alpha)
b(t)\right)+ t\overline{\Pi}^{(-)}(b(t)).
\end{eqnarray*}
Observe that the random variable
$\sum_{0<s\le t}(\Delta_s^{(-)}\wedge b(t))-tA_-(b(t))$ is
centered with
variance $tU_-(b(t))$.
Hence 
$$
P\left(\sum_{0<s\le t}(\Delta_s^{(-)}\wedge
b(t))-tA_-(b(t))>(1/2-\alpha)
b(t)\right)
\leq  \frac{tU_-(b(t))}{(1/2-\alpha)^2b^2(t)},
$$
so that, by the choice of $\ \alpha$, we finally arrive at
\begin{equation}\label{5.4}
P(X_t^{(-)}  > b(t)/2)
\le \left(\frac{2}{(1/2-\alpha)^2}+1\right)\alpha
\le 1/2. \end{equation}

By \eqref{5.2}, 
$P(X_t^{(+)}>b(t)\ {\rm i.o.})\ge P(\Delta_t^{(+)}>b(t)\ {\rm
i.o.})=1$. 
Choose $t_n\downarrow 0$ such that
$P(X_{t_n}^{(+)}>b(t_n)\ {\rm i.o.})=1$. Since the subordinators
$X^{(+)}$
and $X^{(-)}$ are independent, we have
\begin{align*}
&P(X_{t_n}>b(t_n)/2\ {\rm i.o.})
 \nonumber \\
&\ge \lim_{m\to\infty}P(X_{t_n}^{(+)} >b(t_n), X_{t_n}^{(-)}\le
b(t_n)/2
 \ {\rm for\ some}\ n>m)
 \nonumber \\
&\ge (1/2)
P(X_{t_n}^{(+)} >b(t_n)\ {\rm i.o.})\qquad
\ {\rm (by}\ \eqref{5.4} )
\nonumber \\
&=1/2.
 \end{align*}
In the last inequality we used
the Feller--Chung lemma (\cite{CT}, p. 69).
Thus 
$\limsup_{t\downarrow 0} X_t/b(t)\ge 1/2$, a.s.
Now since $a(x)/x$ is nondecreasing, we have for $\alpha<1$,
$b(t)/a(t)\ge a(t/\alpha)/a(t)\ge 1/\alpha$, 
so
$\limsup_{t\downarrow 0} X_t/a(t) \ge 1/\alpha$ a.s.
Letting $\alpha \downarrow 0$ gives 
$\limsup_{t\downarrow 0} X_t/a(t) =\infty$ a.s.,
as claimed in \eqref{a1}.  \hfill \bbox

\medskip

We now state a strong version of the converse of Proposition
\ref{prop3}
which completes the proof of Theorem \ref{thm3.2}.

\begin{prop}\label{prop3-end}
The notation and assumptions are the same as in Proposition
\ref{prop3}.
If \begin{equation}\label{a2-end}
\int_0^1 \frac{\Pi({\rm d}x)}{ 1/a^\leftarrow(x)+A_-(x)/x} <
\infty,\end{equation}
then we have
\begin{equation}\label{a1-end}
\limsup_{t\downarrow 0}
\frac{X_t}{a(t)}\leq 0 \ {\rm a.s.}
\end{equation}
\end{prop}

We will establish Proposition \ref{prop3-end} using
a coupling technique similar to that in \cite{b2}. For this
purpose, we
first need a technical lemma, which is intuitively obvious
once the notation has been assimilated.
Let ${Y}$ be a L\'evy
process
and $((t_i,x_i), i\in I)$ a countable family in
$(0,\infty)\times
(0,\infty)$ such that the $t_i$'s are pairwise distinct. Let
$({Y}^i, i\in I)$ be a family of i.i.d. copies of ${Y}$, 
and set for each $i\in I$
$$\rho_i:=\inf\left\{s\geq0: {Y}^i_s\geq x_i\right\}\wedge
a^\leftarrow(x_i),$$
where $a(\cdot)$ is as in the statement of Proposition
\ref{prop3}.
More generally, we could as well take for $\rho_i$ any stopping
time
in the natural filtration of ${Y}^i$, depending possibly on 
the family $((t_i,x_i), i\in I)$.

 Now assume that 
 \begin{equation}\label{eqcond}
T_t:=\sum_{t_i\leq t}^{}\rho_i\,<\,\infty \hbox{ for all $t\geq
0$ and } \sum_{i\in I}\rho_i=\infty, \quad
\hbox{ a.s.}
\end{equation}
Then $T=(T_t, t\geq0)$ is a right-continuous non-decreasing
process and
(\ref{eqcond}) enables us to construct a process
${Y}'$ by pasting together the paths $({Y}^i_s, 0\leq s \leq \rho_i)$ as
follows.
If $t=T_u$ for some (unique) $u\geq0$, then we set
$${Y}'_t\,=\,\sum_{t_i\leq u} {Y} ^i(\rho_i)\,.$$
Otherwise, there exists a unique $u>0$ such that $T_{u-}\leq
t < T_u$,
and thus  a unique index $j\in I$ for which
$T_u-T_{u-}=\rho_j$,  and  
we set
$${Y}'_t\,=\,\sum_{t_i<u} {Y} ^i(\rho_i)
+{Y}^j(t-T_{u-})\,.$$

\begin{lem}\label{L1} Under the assumptions above, ${Y}'$ is a
version
of ${Y}$; in particular its law does not depend on the family
$((t_i,x_i), i\in I)$.
\end{lem}

\bigskip \noindent{\it Proof of Lemma \ref{L1}:}\
The statement
follows readily from the strong Markov property
in the case when 
the family $(t_i, i\in I)$ is discrete in $[0,\infty)$.
 The general case
is deduced by approximation. \hfill\bbox
\medskip

We will apply Lemma \ref{L1} in the following framework.
Consider a subordinator ${X}^{(-)}$ with no drift and L\'evy
measure $\Pi^{(-)}$; ${X}^{(-)}$ will play the role of the
L\'evy process $Y$ above.
Let also ${X}^{(+)}$ be an independent subordinator with no
drift and L\'evy measure
$\Pi^{(+)}$.
We write $((t_i,x_i), i\in I)$ for the family of the times and
sizes of the jumps  of ${X}^{(+)}$.
 By the L\'evy-It\^o decomposition, $((t_i,x_i), i\in I)$ is the
family
 of the atoms of a Poisson random measure on $\R_+\times \R_+$
with intensity ${\rm d}t\otimes \Pi^{(+)}({\rm d}x)$. 

Next, mark each jump of ${X}^{(+)}$, say $(t_i,x_i)$, using an
independent copy ${X}^{(-,i)}$ of ${X}^{(-)}$.
In other words, 
$\left((t_i,
x_i,{X}^{(-,i)}), i\in I\right)$ is the family of atoms of a
Poisson random measure on $\R_+\times\R_+\times
\D$  with intensity ${\rm d}t\otimes \Pi^{(+)}({\rm d}x)\otimes
\P^{(-)}$,
where $\D$ stands for the space of c\` adl\` ag paths on
$[0,\infty)$ and $\P^{(-)}$ for the law of ${X}^{(-)}$.
Finally, define for every $i\in I$,
$$\rho_i\,:=\,\inf\left\{s\geq0: {X}^{(-,i)}_s\geq
x_i\right\}\wedge
a^\leftarrow(x_i)\,.$$

\begin{lem}\label{L23} In the notation above,
the family $((t_i,\rho_i), i\in I)$ fulfills (\ref{eqcond}).
Further, the process
$$T_t:\,=\,\sum_{t_i\leq t}\rho_i\,,\qquad
t\geq0$$
is a subordinator with no drift.
\end{lem}

\bigskip \noindent{\it Proof of Lemma \ref{L23}:}\
Plainly, 
$$\sum_{i\in I}\delta_{(t_i,
\rho_i)}$$
is a Poisson random measure on $\R_+\times\R_+$ with intensity
${\rm d}t \otimes \mu({\rm d}y)$, where
$$\mu({\rm d}y)\,:=\,\int_{(0,\infty)}\Pi^{(+)}({\rm
d}x)\P^{(-)}\left(\tau_x\wedge
a^\leftarrow(x)\in {\rm d}y\right)\,,$$
and $\tau_x$ denotes the first-passage time of ${X}^{(-)}$ in
$[x,\infty)$.
So it suffices to check that $\int_{(0,\infty)}(1\wedge
y)\mu({\rm d}y)<\infty$.

In this direction, recall (e.g. Proposition III.1 in \cite{b})
that
there is some absolute constant $c$ such that
$$\E^{(-)}\left(\tau_x\right)\leq \frac{cx}{A_-(x)}\,,\qquad
\forall x>0.$$
As a consequence, we have
\begin{eqnarray*}
\int_{(0,\infty)}y\mu({\rm d}y)\,&=&
\,\int_{(0,\infty)}\Pi^{(+)}({\rm d}x)\E^{(-)}\left(\tau_x\wedge
a^\leftarrow(x)\right) \\
&\leq&\,
\int_{(0,\infty)}\Pi^{(+)}({\rm
d}x)\left(\E^{(-)}\left(\tau_x\right)\wedge
a^\leftarrow(x)\right) \\
&\leq&c \int_{(0,\infty)}\Pi^{(+)}({\rm
d}x)\left(\frac{x}{A_-(x)}\wedge
a^\leftarrow(x)\right)\,.
\end{eqnarray*}

Recall that we assume that 
$\Pi^{(+)}$ has support in $[0,1]$.
It is readily checked that convergence of the integral in
\eqref{a2-end} is equivalent to
$$
\int_{(0,\infty)}\Pi^{(+)}({\rm d}x)\left(\frac{x}{A_-(x)}\wedge
a^\leftarrow(x)\right)<\infty.
$$
Our claim is established. \hfill \bbox

We can thus construct a process ${X}'$, as in Lemma \ref{L1},
 by pasting together the paths
$({X}^{(-,i)}_s, 0\leq s \leq \rho_i)$.
This enables us to complete the proof of Proposition
\ref{prop3-end}.

\bigskip \noindent{\it Proof of Proposition \ref{prop3-end}:}\
An application of Lemma \ref{L1}
 shows that ${X}'$ is a
subordinator which is independent of ${X}^{(+)}$ and has the
same
law as
${X}^{(-)}$. As a consequence, we may suppose that the L\'evy
process $X$
is given in the form
$X=X^{(+)}-X'$.

Set $Y_t:={X}'_t+a(t)$. For every jump $(t_i,x_i)$ of
${X}^{(+)}$,
we have by construction
\begin{eqnarray*}
Y(T_{t_i})-Y(T_{t_i-})\,&=&\,{X}^{(-,i)}(\rho_i)+a(T_{t_i-}+\rho_i)-a(T_{t_i-})\\
&\geq&{X}^{(-,i)}(\rho_i)+a(\rho_i)
\qquad \hbox{(as $a(x)/x$\ increases)}\\
&\geq&x_i
\qquad \hbox{(by definition of $\rho_i$)}\,.
\end{eqnarray*}
By summation (recall that ${X}^{(+)}$ has no drift), we get that
$Y(T_{t})\geq {X}^{(+)}_{t}$ for all $t\geq 0$.

As 
$T=(T_t, t\geq0)$ is a subordinator with no drift, we know from
the result
of Shtatland (1965) that $T_t=o(t)$ as $t\to0$, a.s.,
thus  with probability one, we have
for every $\varepsilon>0$
$${X}^{(+)}_t\leq {X}_{\varepsilon t}'+a(\varepsilon t),\
\forall t\geq0
\hbox{ sufficiently small}\,.$$
Since $a(x)/x$\ increases, we deduce that for $t$ sufficiently
small
$$\frac{X_t}{a(t)}\leq \frac{X^{(+)}_t-X'_{\varepsilon
t}}{a(\varepsilon t)/\varepsilon}\leq \varepsilon\,$$
which completes the proof. 
 \hfill\bbox

\subsection{Proof of Theorem  \ref{thm3.3}}

Suppose \eqref{2s} holds, so that $X_t>0$ for all $t\le $ some
(random) $t_0>0$. Thus $X$ is irregular for $(-\infty,0)$ and
\eqref{33b} follows from \cite{b2}.
But \cite{b2} has that \eqref{33b} 
implies
$\sum_{0<s\le t} \Delta^{(-)}_s =o\left(\sum_{0<s\le t}
\Delta^{(+)}_s\right)$, a.s.,
as $t\downarrow 0$, so, for arbitrary $\veps>0$,
$X_t \ge (1-\veps)\sum_{0<s\le t} \Delta^{(+)}_s
:= (1-\veps)X_t^{(+)}$, a.s., when $t\le$
some (random) $t_0(\veps)>0$.
Now $X_t^{(+)}$ is a subordinator with zero drift.
Apply Lemma \ref{ronlem1} with $\gamma=\kappa$ to get
\eqref{33a}.

Conversely, 
\eqref{33b} implies  $\sum_{0<s\le t} \Delta^{(-)}_s
=o\left(\sum_{0<s\le t} \Delta^{(+)}_s\right)$ a.s.,
and \eqref{33a} and Lemma \ref{ronlem1} imply
$\lim_{t\downarrow 0} \sum_{0<s\le t}
\Delta^{(+)}_s/t^\kappa=\infty$ a.s.,
hence  \eqref{2s}.
 \hfill\bbox

\bigskip 
\noindent \textbf{Acknowledgements.}
We are grateful to the participants in the Kiola Research
Workshop, February 2005, where this work began,  for stimulating
input.
The second author also thanks the Leverhulme
Trust for their financial support,
and the second and third authors thank the University of Western
Australia for support and provision of
facilities in the summer of 2006.

\end{document}